\documentclass[11pt, a4paper]{article}

\usepackage{url, hyperref}

\usepackage{amsmath, amssymb, amsthm, mathtools, geometry}

\usepackage{pgfplots}
\usepackage{xfrac}

\theoremstyle{plain}
\newtheorem{theorem}{Theorem}[section]
\newtheorem{conjecture}[theorem]{Conjecture}
\newtheorem{corollary}[theorem]{Corollary}
\newtheorem{lemma}[theorem]{Lemma}
\newtheorem{proposition}[theorem]{Proposition}
\newtheorem*{claim*}{Claim}
\newtheorem*{problem*}{Problem}
\newtheorem*{conjecture*}{Conjecture}

\theoremstyle{definition}
\newtheorem{definition}[theorem]{Definition}
\newtheorem{problem}[theorem]{Problem}
\newtheorem{example}[theorem]{Example}

\newtheorem{question}[theorem]{Question}

\newcommand\al{\alpha}
\newcommand\bt{\beta}

\newcommand\dl{\delta}
\newcommand\lm{\lambda}

\newcommand\la{\langle}
\newcommand\ra{\rangle}
\newcommand\lla{\la\!\la}
\newcommand\rra{\ra\!\ra}

\newcommand\cA{{\cal A}}
\newcommand\cF{{\cal F}}
\newcommand\cJ{{\cal J}}
\newcommand\cM{{\cal M}}

\newcommand\ad{\mathrm{ad}}

\newcommand\FF{\mathbb{F}}
\newcommand\QQ{\mathbb{Q}}
\newcommand\ZZ{\mathbb{Z}}
\newcommand\NN{\mathbb{N}}

\newcommand\Aut{\mathrm{Aut}}

\newcommand{\A}{\mathrm{A}}
\newcommand{\B}{\mathrm{B}}
\newcommand{\C}{\mathrm{C}}
\newcommand{\J}{\mathrm{J}}
\newcommand{\Y}{\mathrm{Y}}
\newcommand{\IY}{\mathrm{IY}}
\renewcommand{\L}{\mathrm{L}}

\newcommand{\Miy}{\mathrm{Miy}}
\newcommand{\ch}{\mathrm{char}}
\newcommand{\Fus}{\mathbf{Fus}}
\newcommand{\Ann}{\mathrm{Ann}}
\newcommand{\supp}{\mathrm{supp}}

\newcommand{\cH}{\mathcal{H}}
\newcommand{\hatH}{\hat{\cH}}

\newcommand{\1}{{\bf 1}}

\renewcommand{\phi}{\varphi}
\renewcommand{\epsilon}{\varepsilon}

\usepackage{enumitem}
\setlist[enumerate,1]{label={\upshape\arabic*.}}
\setlist[enumerate,2]{label={\upshape (\alph*)}}
\setlist[enumerate,3]{label={\upshape (\roman*)}}

\usepackage{xcolor}

\begin{document}

\title{Axial algebras of Jordan and Monster type}

\author{J.~M\textsuperscript{c}Inroy\footnote{Department of Mathematics, University of Chester, Parkgate Rd, Chester, CH1 4BJ, UK, email: J.McInroy@chester.ac.uk}
\and
S.~Shpectorov\footnote{School of Mathematics, University of Birmingham, Edgbaston, Birmingham, B15 2TT, UK, email: S.Shpectorov@bham.ac.uk}}

\date{\today}

\maketitle

\begin{abstract}
Axial algebras are a class of non-associative commutative algebras whose properties are defined in terms of a fusion law.  When this fusion law is graded, the algebra has a naturally associated group of automorphisms and thus axial algebras are inherently related to group theory.  Examples include most Jordan algebras and the Griess algebra for the Monster sporadic simple group.

In this survey, we introduce axial algebras, discuss their structural properties and then concentrate on two specific classes: algebras of Jordan and Monster type, which are rich in examples related to simple groups.
\end{abstract}

\section{Introduction}
Axial algebras are a relatively recent class of non-associative algebras which have a strong natural link to group theory and to physics.  The area has grown significantly since its inception in  2015 and so it is a good time for an article to survey all the recent developments and interesting open problems.

Many interesting classes of algebras, such as associative algebras, Jordan algebras and Lie algebras, are defined by asserting that all the elements should satisfy prescribed polynomial identities.  However, axial algebras are not defined in this way.  Instead, they are defined as being generated by certain idempotents called axes which obey a fusion law.  Obeying an arbitrary generic fusion law here is probably too weak to lead to any interesting theory and the examples that could occur would be wild.  However, specific tight fusion laws, similar to the ones considered in this survey, are both natural enough to allow a rich variety of interesting examples and strong enough to allow them to be classified.

Historically, examples of fusion laws have been known for quite a while.  The Pierce decomposition for associative algebras goes back to the 19\textsuperscript{th} century and later in the 1940s, Pierce decompositions were generalised to Jordan algebras by Albert.  These are in fact nothing other than fusion laws which we will see later in Section \ref{sec:fusion laws}.  However, the fusion language really came to the fore in the 1990s following the discovery of the Monster and the Griess algebra in the 1970s and 80s.  Recall that the Monster is the largest of the 26 sporadic simple groups and it was conjectured to exist by Fischer and Griess in 1973 and shown to exist by Griess in 1982 \cite{G82} as the group of automorphisms of the Griess algebra.  Research around the Moonshine conjecture led to the introduction  of Vertex Operator Algebras (VOAs) by Borcherts \cite{B86} and in particular to the Moonshine VOA $V^\natural$ (the Griess algebra is the weight $2$ part of $V^\natural$).  In 1996, Miyamoto \cite{m} introduced Ising vectors for an OZ-type VOA and he proved that every Ising vector $u$ leads to an involution $\tau_u$ of the VOA.   This is a consequence of the fact that the fusion rules for the representations of Virasoro algebra generated by $u$ are $C_2$-graded.  This idea of associating an involution with an element of the algebra using the fusion law is the start of the axial theory.  (We note that Miyamoto's work built on earlier contributions of Norton and Conway \cite{C85}, who investigated in detail the properties of the Griess algebra.  In particular, they knew of the correspondence between $2\A$ involutions in the Monster and idempotents called $2\A$-axes in the Griess algebra, though without a fusion law.)

Miyamoto began a study of VOAs generated by two Ising vectors and this was completed by Sakuma \cite{s}.  Analysing his proof, Ivanov realised that the calculation is done entirely in the weight $2$-component of the VOA.  He took the properties of this subalgebra used in Sakuma's proof as axioms for a new class of real commutative non-associative algebras called Majorana algebras \cite{i}.  This was the first class where one of the key axioms was a fusion law.  Similarly to the VOA case, special elements in the Majorana algebra, called \emph{Majorana axes}, lead to Majorana involutions in the automorphism group of the algebra.  Initially, the focus was on trying to determine the Majorana algebras for small, particularly simple, groups.  Whilst doing this, it became clear that the fusion law and its grading are the main tools, whilst other axioms are much less important.

Axial algebras were introduced by Hall, Rehren and Shpectorov in \cite{hrs1, hrs2} to further generalise Majorana algebras, removing the unnecessary axioms and allowing arbitrary fusion laws and arbitrary fields.  Historically, most people in this area are group theorists.  This is because of the natural connection to groups tying into the long-standing task of finding a unified theory for simple groups, whereby all of them, including the sporadic groups, could be treated within the same setup.  The class of axial algebras of Monster type, that we discuss below, includes Jordan algebras realising classical groups and $G_2$, Matsuo algebras realising $3$-transposition groups, and of course the Griess algebra and its subalgebras realising the Monster and many other sporadic groups from the Happy Family.  Whilst it might be presumptuous to suppose that axial theory can provide the ultimate unified setup, still we can reasonably expect to find new ties between different kinds of simple groups.

Axial algebra theory has grown significantly since 2015 and we cannot hope to cover everything in this survey.  We are going to introduce axial algebras to the reader and then focus on the structure theory and on two specific classes of axial algebras, those of Jordan and Monster type.

Due to length considerations, we are forced to leave aside some very interesting material.  This includes the link to VOAs and physics, but also recently discovered connections to analysis and other areas of mathematics.  In particular, Tkachev \cite{v} found that algebras arising in the global geometry and regularity theory of non-linear PDEs are axial algebras for suitable compact fusion laws, and similarly Fox \cite{fox} shows axial properties of the algebra of curvature tensors.  Another recent development that we do not cover is the theory of decomposition algebras \cite{dpsv} due to De Medts, Peacock, Shpectorov and Van Couwenberghe, which generalises the axial setup even further and provides a categorical point of view.  Rowen and Segev in \cite{non comm} suggest a partial generalisation to the non-commutative case. Even though Majorana algebras are squarely within the class of axial algebras of Monster type, the motivation and approach of Majorana theory is slightly different and most importantly Ivanov himself published two recent surveys on this topic \cite{Majorana future, Majorana future2}.   Finally, it is very likely we missed some other interesting work and if so, we would like to apologise for this.

\section{Background}
We begin by introducing fusion laws and axial algebras and giving some examples.  We then explain how we get a naturally associated group of automorphisms called the Miyamoto group.

\noindent{\bf Notation:} For an algebra $A$ and a set of elements $X \subseteq A$, we will write $\la X \ra$ for the subspace of $A$ spanned by $X$ and $\lla X \rra$ for the subalgebra generated by $X$.

\subsection{Fusion laws}\label{sec:fusion laws}

The starting point is the following concept.

\begin{definition} \label{fusion law}
A \emph{fusion law} is a pair $\cF := (\cF, \star)$, where $\cF$ is a set and $\star$ is a map, 
\[
\star\colon \cF \times \cF \rightarrow 2^{\cF},
\]
where $2^\cF$ is our notation for the set of all subsets of $\cF$.
\end{definition}
 
The fusion laws we are most interested in have a small set $\cF$, and it will be convenient to present them in a table similar to a group multiplication table. For example, in Table \ref{tab: fusion laws}, we show the three fusion laws: the associative law $\mathcal{A}$, the Jordan type law $\cJ(\eta)$ and the (generalised) Monster type law $\cM(\al, \bt)$.  For simplicity, we drop the set signs and simply list the elements of the set $\lm \star \mu$ in the $(\lm, \mu)$ cell.  In particular, the empty cell represents the empty set. 

\begin{table}[h]

\setlength{\tabcolsep}{4pt}
\renewcommand{\arraystretch}{1.6}
\centering
	\begin{minipage}[t]{0.2\linewidth}
		\begin{tabular}{c||c|c}
		$\star$ & $1$ &$0$\\
		\hline \hline
		$1$&$1$& \\
		\hline
		$0$& &$0$
		\end{tabular}
	\end{minipage}
	\begin{minipage}[t]{0.25\linewidth}
		\begin{tabular}{c||c|c|c}
		$\star$ & $1$ &$0$&$\eta$\\
		\hline  \hline
		$1$&$1$& &$\eta$\\
		\hline
		$0$& &$0$&$\eta$\\
		\hline
		$\eta$&$\eta$&$\eta$&$1,0$
		\end{tabular}
	\end{minipage}
	\begin{minipage}[t]{0.25\linewidth}
		\begin{tabular}{c||c|c|c|c}
		$\star$ & $1$ &$0$&$\alpha$& $\beta$\\
		\hline \hline
		$1$&$1$& &$\alpha$& $\beta$\\
		\hline
		$0$& &$0$&$\alpha$& $\beta$\\
		\hline
		$\alpha$&$\alpha$&$\alpha$&$1,0$& $\beta$\\
		\hline
		$\beta$&$\beta$&$\beta$&$\beta$&$1,0, \alpha$
		\end{tabular}
	\end{minipage}
\caption{The $\cA$, $\cJ(\eta)$, $\cM(\alpha,\beta)$ fusion laws.}
\label{tab: fusion laws}
\end{table}

Note that these examples are all \emph{symmetric} fusion laws, i.e.\ $\lm \star \mu = \mu \star \lm$, for all $\lm, \mu \in \cF$.  From now on, we will consider only symmetric fusion laws.

\subsection{Axial algebras}

Let $A$ be a commutative algebra over a field $\FF$. An adjoint map $ad_a$ of an element $a\in A$
is the endomorphism of the algebra which takes $x$ to $ax$. For $\lambda \in \FF$, 
\[
A_{\lambda}(a)=\{ x \in A : ax=\lambda x \}
\]
is the corresponding eigenspace. Note that $A_{\lambda}(a)=0$ if $\lambda$ is not an eigenvalue of the adjoint map $ad_a$.  We further extend the eigenspace notation to sets $S \subseteq \FF$ by defining $A_S(a) := \bigoplus_{\lm \in S} A_\lm(a)$.

We use fusion laws to impose partial control on the multiplication in an algebra.

\begin{definition}
Let $A$ be a commutative, nonassociate algebra over a field $\FF$ and $(\cF, \star)$ be a fusion law with $\cF \subseteq \FF$.
An element $a$ is said to be an \emph{$\cF$-axis} (or just an axis if $\cF$ is clear from context) if
\begin{enumerate}
\item
$a$ is an idempotent, that is $a^2=a$;
\item
the adjoint map $ad_a$ is semisimple with all eigenvalues contained in $\cF$, that is,
\[
A= A_\cF(a);
\]
\item
the above decomposition obeys the fusion law $\cF$.  That is, for all $\lambda, \mu \in \cF$,
\[
A_{\lambda}(a) A_{\mu}(a) \subseteq A_{\lambda \star \mu}(a).
\]

\end{enumerate}
\end{definition}

From the above definition, as $a$ is an idempotent, $1$ is an eigenvalue of $ad_a$.  So we always assume that $1 \in \cF$.  

\begin{definition}
An axis $a \in A$ is said to be \emph{primitive} if $A_1(a)=\FF a$.
\end{definition}

\begin{definition}
A \emph{primitive $\cF$-axial algebra} is a pair $A = (A, X)$, where $A$ is a commutative, non-associative algebra over the field $\FF$ and $X$ is a set of primitive $\cF$-axes which generate $A$.
\end{definition}

We will (mostly) just consider primitive axes and axial algebras and so for simplicity we will typically just speak of axes and axial algebras assuming primitivity.

It is also important to note that $X$ is just \textit{a} set of axes.  We never assume that $X$ is the set of all idempotents that have the property of axes.  In fact, the problem of finding all such idempotents in a known algebra $A$ is a very difficult problem.  This question is related to finding the full automorphism group of an axial algebra, which we discuss later in Section \ref{sec:aut}.

We use the fusion law as a property which defines a class of axial algebras.  For example, $\cJ(\eta)$- axial algebras constitute the class of axial algebras of Jordan type $\eta$ and $\cM(\al, \bt)$-axial algebras constitute the class of algebras of Monster type $(\al,\bt)$.  The following example was one of the motivations for defining axial algebras.

\begin{example}
The Monster sporadic simple group $M$ was constructed as the automorphism group of a $196,884$-dimensional real algebra $V$, known as the Griess algebra.  The algebra $V$ together with the set of so-called $2\A$-axes in it is an axial algebra of Monster type $(\frac{1}{4}, \frac{1}{32})$.
\end{example}

Recall that a Jordan algebra is a commutative non-associative algebra $A$ satisfying the \emph{Jordan identity}
\[
(xy)x^2 = x(y x^2)
\]
for all $x,y \in A$.  The class of Jordan algebras also gives us examples of axial behaviour, namely every idempotent (primitive or not) in a Jordan algebra is a $\cJ(\frac{1}{2})$-axis.  We should be careful here, because not every Jordan algebra is generated by non-zero idempotents.  However, if a Jordan algebra $A$ is generated by primitive idempotents, then it is an axial algebra of Jordan type $\frac{1}{2}$.  In particular, all simple Jordan algebras over an algebraically closed field are examples.

One example which we will see again later is the spin factor Jordan algebra, which can be defined as follows.

\begin{example}\label{spin factor}
Let $V$ be a vector space over a field $\FF$ of characteristic not $2$ and $b \colon V \times V \to \FF$ be a symmetric bilinear form.  Let $A = S(b) = \FF \1 \oplus V$ and define multiplication on $A$ by
\[
\begin{gathered}
\1^2 = \1, \qquad \1u = u \\
uv = \tfrac{1}{2}b(u,v) \1
\end{gathered}
\]
for all $u,v \in V$.  This algebra is known as the \emph{spin factor} algebra, or the Jordan algebra of \emph{Clifford type}.  The (primitive) axes have the form $\frac{1}{2}(1+u)$, where $b(u,u) = 2$, and so the spin factor is an axial algebra of Jordan type $\frac{1}{2}$ if and only if $V$ is spanned by vectors $u$ with $b(u,u) = 2$.
\end{example}

Here is another interesting example of algebras of Jordan type that we will see again later.

\begin{example}\label{Matsuo algebra}
Let $(G, D)$ be a group of $3$-transpositions.  That is, $D$ is a normal set of involutions in $G$ such that $G = \la D \ra$ and $|cd| \leq 3$ for all $c,d \in D$.  For example, we can take $G = S_n$ and $D$ the set of all transpositions in $G$.

Suppose $\FF$ is a field of characteristic not $2$ and $\eta \in \FF - \{0,1\}$.  The \emph{Matsuo algebra} $A := M_\eta(G,D)$ has basis $D$ and multiplication given by
\[
a \circ b = \begin{cases}
a & \mbox{if } b= a, \\
0 & \mbox{if } |ab| =2, \\
\frac{\eta}{2}(a+b-c) & \mbox{if }  |ab| =3, \mbox{ where } c = a^b = b^a,
\end{cases}
\]
where we use $\circ$ for the algebra product to distinguish it from the multiplication in $G$.  It is easy to see that for $a \in D$, $A_1(a) = \la a \ra$, $A_0(a) = \la b : |ab| = 2 \ra \oplus \la b+ c - \eta a : |ab| = 3 \ra$ and $A_\eta(a) = \la b- c : |ab| = 3 \ra$, where as above $c = a^b = b^a$.  By counting dimensions, we can see that $\ad_a$ is semisimple.  A short calculation shows that the eigenspaces obey the $\cJ(\eta)$ fusion law and so every $a \in D$ is a (primitive) $\cJ(\eta)$-axis.  Since the set of all $a \in D$ generate the algebra $A$, it is an axial algebra of type $\cJ(\eta)$.
\end{example}

Let us finish this subsection with a brief discussion of what additional assumptions we can make about the fusion laws due to our axial algebra definitions.

First of all, we have already mentioned that we always assume that $1\in\cF$, because axes are idempotents. Furthermore, our assumption that the fusion law be symmetric is natural because of commutativity of axial algebra. If we were to consider non-commutative generalisations, we would also have to allow non-symmetric fusion laws. Finally, primitivity also has some implications. Namely, 
for a primitive axis $a\in A$, we have that $A_1(a)=\la a\ra$ and, because of this, 
\[
A_1(a)A_\lm(a)=aA_\lm(a)=\begin{cases}A_\lm(a)&\mbox{if }\lm\neq 0;\\0&\mbox{if }\lm=0.\end{cases}
\]
So, for primitive axial algebras, we always assume that $1\star\lm=\{\lm\}$ if $\lm\neq 0$ and $1\star 0=\emptyset$.


\subsection{Miyamoto group}\label{sec:Miy}

The important feature of axial algebras is that they are closely related to groups. One can see this in the examples we presented earlier: The Griess algebra is clearly related to the Monster group, which is its automorphism group. The Matsuo algebras are clearly related to $3$-transposition groups. Also, the Jordan algebras are very symmetric allowing classical and even some exceptional groups of Lie type as their automorphism groups. This is not an accident because the fusion laws of these axial algebras (see Table \ref{tab: fusion laws}) are graded by the group $C_2$. 

The exact meaning of this is as follows. For the law $\cF=\cJ(\eta)=\{1,0,\eta\}$, let the plus part of it be $\cF_+=\{1,0\}$ and the minus part be $\cF_-=\{\eta\}$. Then we have from the fusion law that $\cF_+\star\cF_+ \subseteq\cF_+$, $\cF_-\star\cF_-\subseteq\cF_+$, and $\cF_+\star\cF_-\subseteq\cF_-$. 
In an algebra $A$ of Jordan type $\eta$, this results in a $C_2$-grading for every axis $a$: setting $A_+:=A_{\cF_+}(a)$ and $A_-:=A_{\cF_-}(a)$,  we have that
\[
A_+A_+\subseteq A_+, A_-A_-\subseteq A_+,\mbox{ and }A_+A_-\subseteq A_-.
\]
It is well known that such a $C_2$-grading of the algebra leads to an automorphism $\tau_a$ of $A$. This $\tau_a$ acts as identity on $A_+$ and as minus identity on $A_-$. (Then, $\tau_a^2 = 1$ and $\tau_a = 1$ only if $A_-=0$, which is an exceptional situation.) We stress that we get a separate  involution for each axis $a$, and thus, since $A$ typically contains many axes, this leads to a substantial group of automorphisms.

For example, in the case of Matsuo algebras, this construction allows us to recover the $3$-transposition group the algebra was constructed from, albeit we only recover the group up to its centre. In the case of Jordan algebras, we recover in this way a class of involutions in the corresponding group of Lie type.

The same grading trick works for the fusion law $\cF=\cM(\al,\bt)$, where we split $\cF$ as follows: $\cF_+=\{1,0,\al\}$ and $\cF_-=\{\bt\}$. Again, this is a $C_2$-grading of the fusion law and it results in $C_2$-gradings and corresponding involutions for all axial algebras of Monster type. In the Griess algebra, these involutions are the $2\A$ involutions in the Monster and thus we recover the Monster $M$ in a natural way from its algebra.

Let us now give the formal definitions. Initially, for example, in \cite{hrs1,hrs2}, the focus was specifically on the case of a $C_2$-grading, as above, but a general abelian grading group $T$ can be handled in pretty much the same way, in terms of partitions. For example, such a definition was used in \cite{s}. (A more categorical approach to gradings will be discussed later.)

\begin{definition}
Suppose $\cF$ is a fusion law and $T$ is an abelian group. A $T$-grading of $\cF$ is a partition $\{\cF_t: t\in T\}$ (with some parts $\cF_t$ possibly empty) such that $\cF_s\star\cF_t\subseteq \cF_{st}$ for all $s,t\in T$.
\end{definition}

Given a $T$-grading of $\cF$ and an axis $a$ in an $\cF$-axial algebra $A$, it is immediate that 
\[
A=\bigoplus_{t\in T}A_{\cF_t}(a)
\] 
is a $T$-grading of the algebra $A$. It now follows that, for each linear character $\chi$ of $T$, we can define an automorphism of $A$ as follows. Let $\tau_a(\chi)$ be the linear map $A\to A$, which acts on $A_t(a):=A_{\cF_t}(a)$ by multiplying with the scalar $\chi(t)$.

The automorphism $\tau_a(\chi)$ is called a \emph{Miyamoto automorphism}. It is immediate that the map $\chi\mapsto\tau_a(\chi)$ is a homomorphism from the group $T^\ast$ of linear characters of $T$ to $\Aut(A)$. The image of this map, $T_a\leq\Aut(A)$ is called the \emph{axial subgroup}.

\begin{definition}
For a $T$-graded fusion law $\cF$ and an $\cF$-axial algebra $A$ with the set of generating axes $X$, the group 
\[
\Miy(A, X)=\la T_a:a\in X\ra=\la\tau_a(\chi):a\in X,\chi\in T^\ast\ra\leq\Aut(A)
\]
is called the \emph{Miyamoto group} of $A$.  We will also use $\Miy(A)$ and $\Miy(X)$, depending on what our focus is.
\end{definition}

In the case of $T=C_2$, when the characteristic of the ground field $\FF$ is not $2$, the group $T^\ast$ contains a unique non-trivial character $\chi$ and so in this case we simplify $\tau_a(\chi)$ to simply $\tau_a$ and call it the Miyamoto involution, as above.

The paper \cite{dpsv} introduced the category of fusion laws allowing us to separate the fusion law from the ground field $\FF$. In this, more modern language, a grading is simply a morphism from the fusion law $\cF$ to the group fusion law $T$. Let us provide the exact definitions.

The category $\Fus$ has as objects all fusion laws, as given in Definition \ref{fusion law}, and morphisms are as follows.  

\begin{definition}
Let $(\cF_1, \star_1)$ and $(\cF_2, \star_2)$ be fusion laws.  A \emph{morphism} from $(\cF_1, \star_1)$ to $(\cF_2, \star_2)$ is a map $f \colon \cF_1 \to \cF_2$ such that
\[
f(\lm \star_1 \mu) \subseteq f(\lm) \star_2 f(\mu)
\]
for all $\lm, \mu \in \cF_1$.
\end{definition}

The category of groups is a full subcategory of $\Fus$ via the following construction.

\begin{definition}
For a group $T$, we define the \emph{group fusion law} $(T,\star)$, where $s\star t=\{st\}$ for all $s,t\in T$. 
\end{definition}

It is easy to see that morphisms of group fusion laws are the same as group homomorphisms, so here we indeed have a full subcategory of $\Fus$.

In this categorical language, a $T$-grading of a fusion law $\cF$ becomes simply a morphism from $\cF$ to the group fusion law of $T$. (The parts $\cF_t$ are then the fibres of the morphism.) In this sense, we don't even need to assume that $T$ is abelian. However, since our algebras are commutative and Miyamoto automorphisms come from linear characters, we cannot gainfully utilise such a full generalisation.

Note that for a $T$-graded fusion law $\cF$, we can always unnecessarily enlarge the group $T$ and still have a grading. We say the grading $f\colon\cF\to T$ is \emph{adequate} if $f(\cF)$ generates $T$. Clearly, we can restrict ourselves to only adequate gradings.  A grading $f\colon\cF\to T$ is a \emph{finest} grading if every other grading of $\cF$ factors uniquely through $f$.  Every fusion law admits a unique finest grading and this is adequate \cite{dpsv}.

Notice that if $g \in \Aut(A)$ and $a$ is an axis, then $a^g$ also has all the properties require for an axis and it is primitive if and only if $a$ is primitive.  So we can add it to our set $X$ of generating axes.  If we do this with $g \in \Miy(X)$, this gives us the concept of the \emph{closure} $\bar{X} = X^{\Miy(X)} = \{ x^g : g \in \Miy(X), x \in X \}$ of the set $X$ of generators.

\begin{lemma}
The Miyamoto group $\Miy(\bar{X}) = \Miy(X)$.  In particular, $\bar{\bar{X}} = \bar{X}$.
\end{lemma}

So our notation $\Miy(A)$ is unambigous even though we allow an extension of the generating set $X$ to its closure as above.  Furthermore, in what follows, we will often assume that $X$ is \emph{closed} i.e.\ that $X = \bar{X}$.  This has the advantage that we have an action of $\Miy(A)$ on the set $X$ of generating axes.  This action is necessarily faithful as the axes generate the algebra and so any automorphism fixing all the axes acts trivially on the algebra.

\subsection{Frobenius form}

Often algebras admit bilinear forms that associate with the algebra product. For example, the Killing form on a Lie algebra has this property.

\begin{definition}
Suppose $A$ is an axial algebra. A non-zero bilinear form $(\cdot,\cdot)$ on $A$ is called a \emph{Frobenius} form if 
\[
(u,vw)=(uv,w)
\]
for all $u,v,w\in A$.
\end{definition}

This terminology came from the theory of Frobenius algebras which are associative algebras admitting a non-degenerate form satisfying the above associativity condition. Note that we don't assume that a Frobenius form is necessarily non-degenerate. In fact, as we will see later, the radical of the Frobenius form plays an important role in the structure theory of axial algebras.

Also note that in some areas of mathematics, where algebras with axial properties arise, algebras admitting a Frobenius form are called metrizable.

Often a Frobenius form would have additional nice properties. For example, the Griess algebra admits a Frobenius form which is positive definite and each axis has norm $1$. These properties became some of the axioms of Majorana algebras. Similarly, every Jordan algebra admits a Frobenius form with respect to which every primitive idempotent has norm $1$.

Let us mention a couple of general facts about Frobenius forms on axial algebras.

\begin{proposition} \label{basic Frobenius}
\begin{enumerate}
\item Every Frobenius form is symmetric.
\item For an axis $a\in A$, the eigenspaces $A_\lm(a)$ and $A_\mu(a)$, where $\lm\neq\mu$, are orthogonal with respect to every Frobenius form.
\end{enumerate}
\end{proposition}

As of today, we actually don't know any examples of axial algebras which don't admit a Frobenius form. It is hard to imagine that no such algebras exist.

\section{Structure theory}
The structure theory of axial algebras was developed in \cite{kms1} and the results in this section are from there unless otherwise referenced.

\subsection{Radical}

\begin{theorem} \label{unique ideal}
Suppose that $(A,X)$ is a primitive axial algebra. Then there exists a unique maximal ideal disjoint from the set $X$ of generating axes.
\end{theorem}

\begin{definition} \label{radical}
The unique maximal ideal from Theorem \ref{unique ideal} is called the \emph{radical} of the axial algebra $A$.
\end{definition}

Our notation for the radical is $R(A,X)$. Note that the radical a priori depends on the set of generating axes. However, in most cases, this set $X$ is assumed, and then we just write $R(A)$.

The above definition is implicit, which is an inconvenience. When the algebra admits a Frobenius form, this can be remedied. 

\begin{theorem} \label{when Frobenius}
Suppose that $A$ admits a Frobenius form such that $(a,a)\neq 0$ for all $a\in X$. Then $R(A,X)$ coincides with the radical $A^\perp$ of the Frobenius form.
\end{theorem}

Note that if $(a,a)=0$ for some $a\in X$ then, as $a$ is primitive, it follows from Proposition \ref{basic Frobenius} that $a\perp A$, that is, $a$ is contained in the radical of the form. So the condition in the theorem is a necessary one.

Clearly, the above theorem gives us a practical and easy way to compute the radical by means of linear algebra. When the form is non-degenerate, the radical is trivial. In particular, this is the case for Majorana algebras where the Frobenius form is assumed to be positive definite. Note also that the condition on axes from the above theorem becomes mute for a positive definite form, so in a Majorana algebra every non-zero ideal contains a generating primitive axis.

In fact, this observation can be generalised.

\begin{theorem}
If $(A,X)$ is a primitive axial algebra admitting a non-degenerate Frobenius form then $R(A,X)=0$, i.e., every non-zero ideal contains a generating primitive axis.
\end{theorem}

Indeed, as we already pointed out, if $(a,a)=0$ for a primitive axis $a\in A$, then $a\in A^\perp$. Since $A^\perp=0$, we have that the assumptions of Theorem \ref{when Frobenius} are satisfied and so $R(A,X)=A^\perp=0$.

\subsection{Connectivity}

In the previous subsection, we looked at ideals which do not contain any axes, so now we consider ideals which do contain axes.  We begin with some observations which apply to any ideal and in fact are used in the proofs of the results from the previous subsection.

Since an ideal $I$ is invariant under $\ad_a$ for all axes $a$ and $\ad_a$ is semisimple, then $\ad_a$ is also semisimple acting on $I$ and we have
\[
I = \bigoplus_{\lm \in \cF} I_\lm(a)
\]
where $I_\lm(a) := A_\lm(a) \cap I$.  In particular, this means that $I$ is invariant under every Miyamoto automorphism $\tau_a(\chi)$ and hence every ideal $I$ is invariant under the Miyamoto group $\Miy(A)$.  This will be generalised later.

Suppose that $a$ is an axis and $v \in A$.  We can decompose $v$ with respect to the eigenspace decomposition given by $\ad_a$
\[
v = \sum_{\lm \in \cF} v_\lm
\]
where $v_\lm \in A_\lm(a)$.  In particular, if $a$ is a primitive axis, then $v_1 = \phi_a(v) a$ for some $\phi_a(v) \in \FF$.  We call $\phi_a(v)$ the \emph{projection} of $v$ onto $a$.  

In particular, if $v$ is in an ideal $I$, then $v_\lm \in I$ for all $\lm \in \cF$.  If furthermore $a$ is primitive, then we have the following dichotomy: either the projection $\phi_a(v) = 0$ and so $v_1 = 0$, or the projection $\phi_a(v) \neq 0$ and then $ 0 \neq v_1 = \phi_a(v) a \in I$, meaning $a \in I$.  We will now apply this dichotomy with $v = b$ an axis in $A$ to get the following definition.

\begin{definition}
The \emph{projection graph} $\Gamma$ is a directed graph whose vertices are the axes $X$ and there is a directed edge $b \leftarrow a$ if $\phi_a(b) \neq 0$.
\end{definition}

This immediately gives us the following lemma.

\begin{lemma}\label{proj graph}
If an axis $a \in I$ and there exists a directed path from $a$ to $b$ in the projection graph $\Gamma$, then $b \in I$.
\end{lemma}

In other words, if $a$ is in $I$, then all the successors of $a$ in the projection graph are also in $I$.  If $\Gamma$ is strongly connected, i.e.\ every vertex is a successor to every other vertex, then $A$ has no proper ideals containing an axis.

Assuming the set $X$ of axes is closed under the action of the Miyamoto group, then we can form the \emph{orbit projection graph} $\bar{\Gamma} := \Gamma / \Miy(A)$ as the quotient of the projection graph $\Gamma$.  We then have the analogous result in Lemma \ref{proj graph} for the orbit projection graph $\bar{\Gamma}$.

If we our algebra has a Frobenius form, then the picture is further simplified.

\begin{lemma}
Let $A$ be a primitive axial algebra with a Frobenius form and suppose that $(a,a) \neq 0 \neq (b,b)$, for axes $a$ and $b$.  Then the following are equivalent.
\begin{enumerate}
\item There is a directed edge $a \rightarrow b$.
\item There is a directed edge $a \leftarrow b$.
\item $(a,b) \neq 0$.
\end{enumerate}
\end{lemma}

So in particular, if the Frobenius form is non-zero on all the axes, then we can consider the projection and orbit projection graphs as undirected graphs.

One further observation about ideals which contain axes is the following.  If $a \in I$, then it is easy to see that $A_\lm(a) \subseteq I$ for all $\lm \neq 0$.

\subsection{Associativity}

Although axial algebras usually have elements $x,y,z \in A$ where $(xy)z \neq x(yz)$, they can be associative algebras (so non-associative for us means not necessarily associative).  In fact, it is easy to characterise the associative axial algebras.

\begin{theorem}\textup{\cite[Corollary 2.9]{hrs2}}
Let $A$ be a \textup{(}primitive\textup{)} $\cA$-axial algebra $A$ generating by a set of axes $X$.  Then $A$ is associative and $A = \FF X$.  Conversely, if $A = \lla X \rra$ is an associative axial algebra, then every axis $a \in X$ satisfies the $\cA$ fusion law.
\end{theorem}

This justifies our calling the $\cA$ fusion law (see Table \ref{tab: fusion laws}) the associative fusion law.

Even if our axial algebra is not associative, if the fusion law has the so-called Seress property, then the axes have some measure of commutativity.

\begin{definition}
A fusion law $\cF$ is \emph{Seress} if $0 \in \cF$ and for all $\lm \in \cF$, $0 \star \lambda \subseteq \{ \lambda\}$.
\end{definition}

Since axes are idempotents and the $\lambda \in \cF$ are eigenvalues for the adjoint $\ad_a$, $1 \in \cF$ and $1 \star \lambda \subseteq \{ \lambda \}$.  So for a Seress fusion law, we have $1 \star 0 \subseteq \{0\} \cap \{1\} = \emptyset$.  Also, for a Seress fusion law $0 \star 0 \subseteq \{0 \}$ and so $A_0(a)$ is a subalgebra of $A$ for any axis $a \in X$.  Looking at Table \ref{tab: fusion laws} of examples of fusion laws, we can see that they are all Seress.

\begin{lemma}[Seress Lemma]\footnote{This was first proved by Seress for axial algebras of Monster type $(\frac{1}{4}, \frac{1}{32})$, but Hall, Rehren and Shpectorov \cite[Proposition 3.9]{hrs1} noticed that the same proof held in general.}
Let $A$ be an axial algebra with a Seress fusion law and $a \in X$.  Then for all $y \in A_1(a) \cup A_0(a)$ and $x \in A$, we have
\[
a(xy) = (ax)y
\]
\end{lemma}

\subsection{Sum decompositions}

The subject of axes and ideals naturally leads onto direct sums and more generally sum decompositions of axial algebras.  

\begin{definition}
An algebra $A$ has a \emph{sum decomposition} if it is generated by a set $\{ A_i : i \in I \}$ of pairwise annihilating subalgebras i.e.\ $A_i A_j = 0$ for all $i \neq j$.  We write $A = \Box_{i \in I} A_i$.
\end{definition}

Note that we immediately have that $A = \sum_{i \in I} A_i$ as a vector space decomposition if $A$ has a sum decomposition.  If $A$ has a sum decomposition and this sum is direct, then we write $A = \boxplus_{i \in I} A_i$.

Recall that the \emph{annihilator ideal} is $\Ann(A) := \{ x \in A : xA = 0 \}$.  For an axial algebra, as axes are idempotents, $\Ann(A)$ does not contain any axes and so $\Ann(A) \subseteq R(A)$.

\begin{proposition}
Let $A = \Box_{i \in I} A_i$ be an axial algebra.
\begin{enumerate}
\item $\Ann(A) = \Box_{i \in I} \Ann(A_i)$.
\item If $\Ann(A) = 0$, then $A = \boxplus_{i \in I} A_i$.
\end{enumerate}
\end{proposition}

So for the sum not to be direct, at least one of the $A_i$ must have non-trivial annihilator.  In fact, it can be shown that at least two of the $A_i$ must have non-trivial annihilator.

Just like for (not necessarily commutative) algebras we define the \emph{commutator} $[x,y] := xy - yx$ to measure commutativity, for a non-associative algebra $A$ we define the \emph{associator} $(x,y,z) := (xy)z - x(yz)$ (see for example \cite{Schafer}).  The \emph{centre} $Z(A)$ is the subalgebra
\[
Z(A) := \{ a \in A : [a,A] = 0, (a,A,A) = (A,a,A) = (A,A,a) = 0 \}
\]
In fact, it is easy to see that if the algebra is commutative (as axial algebras are), then $(a,A,A) = 0$ is equivalent to $(A,A,a)=0$ and either of these imply $(A,a,A)=0$.

\begin{proposition}[Turner\footnote{Turner is currently a PhD student at the University of Birmingham and this result is part of his MRes (Qual) project.}, 2022]
Let $A$ be an axial algebra.
\begin{enumerate}
\item For every axis $a$, $Z(A) \subseteq A_1(a) \oplus A_0(a)$.
\item If $A$ is primitive, then $\Ann(A) = Z(A) \cap R(A,X)$.
\end{enumerate}
\end{proposition}

Suppose that we have an axial algebra $A$ which has a sum decomposition $A = \Box_{i \in I} A_i$.  We have so far not put any restriction on the $A_i$ -- these subalgebras might not be axial subalgebras.  That is, they might not be axial algebras in their own right.

\begin{theorem}
Suppose that $A = \Box_{i \in I} A_i$ is an axial algebra generated by a set of axes $X$.  Let $X_i := X \cap A_i$ and define $B_i := \lla X_i \rra$.  Then the $X_i$ partition the set $X$ and $A = \Box_{i \in I} B_i$.
\end{theorem}

So whenever we have a sum decomposition of an axial algebra, we in fact have a sum decomposition of $A$ with respect to axial subalgebras of $A$.  Note however, that we may have $A_i \neq B_i$.

It is clear that if we have two different sum decomposition of an axial algebra into axial subalgebras, we can take a refinement of the two simply by taking a refinement of the two partitions of the set $X$.  This raises the question of what the finest sum decomposition of an axial algebra is.  In order to tackle this, we introduce the following graph.

\begin{definition}
The \emph{non-annihilating} graph $\Delta(X)$ of an axial algebra generated by a set of axes $X$ has $X$ as its vertex set and an edge $a \sim b$ between distinct vertices if $ab \neq 0$.
\end{definition}

This was first introduced for axial algebras of Jordan type $\eta$ in \cite{hss1} and for arbitrary axial algebras in \cite{kms1}.  In fact, if $A$ is an axial algebra which has a Frobenius form which is non-zero on all the axes, then the non-annihilating graph is just the same as the projection graph.

If $A$ has a sum decomposition $\Box_{i \in I} A_i$, then by definition $A_i A_j = 0$ and hence $X_i X_j = 0$.  So each $X_i$ must be a union of connected components in the non-annihilating graph $\Delta(X)$.  This naturally given us the following conjecture.

\begin{conjecture}\label{conj:sum decomp}
The finest sum decomposition of an axial algebra $A$ arises when the $X_i$ are the connected components of the non-annihilating graph $\Delta(X)$.
\end{conjecture}

In other words, let the $X_i$ be the connected components of the non-annihilating graph and set $B_i := \lla X_i \rra$.  It is clear that the $B_i$ generate $A$ as $A$ is generated by the set $X = \bigcup_{i \in I} X_i$.  So to prove the conjecture, we need to show that $B_iB_j = 0$ for all $i \neq j$ and hence $A = \Box_{i \in I} B_i$.

This conjecture holds for axial algebras of Jordan type $\eta$ \cite[Theorem A]{hss1}.  Before we say what progress has been made for this conjecture, note that if $A$ has a non-trivial sum decomposition with axes in different parts, then $0 \in \cF$.  Likewise, if $\Delta(X)$ is not connected, then $0 \in \cF$.

Suppose that the fusion law is graded and so the axial algebra has a non-trivial Miyamoto group.  This group must behave well with respect to the finest sum decomposition.

\begin{theorem}
Suppose $A$ is a $T$-graded axial algebra, with $0 \in \cF_{1_T}$.  Let $X_i$, $i \in I$, be the connected components of $\Delta(X)$ and set $B_i := \lla X_i \rra$.  Then, $\Miy(A)$ is a central product of the $\Miy(B_i)$.
\end{theorem}

Before stating our partial result for Conjecture \ref{conj:sum decomp}, we must first introduce quasi-ideals which generalise the notion of ideals.

\begin{definition}
Let $A$ be an axial algebra with set of generating axes $X$.  A \emph{quasi-ideal} is a subspace $I \subseteq A$ such that $aI \subseteq I$ for all $a \in X$.
\end{definition}

\begin{definition}
The \emph{spine} of an axial algebra $A$ is the quasi-ideal $Q(A, X)$ generated by the axes $X$.  We say $A$ is \emph{slender} if $Q(A,X)=A$.
\end{definition}

We can now state our partial result for the conjecture.

\begin{theorem}
Suppose $A$ is an axial algebra with a Seress fusion law.  Let $X_i$ be the connected components of $\Delta(X)$ and set $A_i := \lla X_i \rra$.  If all but possibly one of the $A_i$ are slender, then $A = \Box_{i \in I} A_i$ and so the conjecture holds.
\end{theorem}

In fact, we know of very few examples of axial algebras which are not slender.

\section{Axial algebras of Jordan type}
An axial algebra is of \emph{Jordan type $\eta$}, or just of Jordan type, if it has the Jordan fusion law $\cJ(\eta)$.  This class was introduced in 2015 in \cite{hrs1} and is one of the most widely studied.  We saw in Examples \ref{spin factor} and \ref{Matsuo algebra} several examples of axial algebras of Jordan type.  One of these was related to $3$-transposition groups.

\subsection{$3$-transposition groups}

Recall that a $3$-transposition group is a pair $(G,D)$ of a group $G$ generated by a normal set of involutions $D$ such that $|cd| \leq 3$ for all $c,d \in D$.  The \emph{diagram} on a set $X\subseteq D$ is the graph having $X$ as a vertex set, where $a,b\in X$ are adjacent whenever $|ab|=3$ (so the complement of the diagram is the commuting graph). We say that $(G,D)$ is \emph{connected} when the diagram on $D$ is connected.\footnote{Hall uses the term ``$3$-transposition group'' only in the connected case.} It is not hard to see that $(G,D)$ is connected if and only if $D$ is a single conjugacy class.  A general $3$-transposition group decomposes as a central product of connected $3$-transposition groups corresponding to the conjugacy classes within $D$, so it is enough to classify the connected $3$-transposition groups.

Fischer began this, when he classified the connected irreducible $3$-transposition groups.  Here irreducibility is a technical condition that requires $[G,O_2(G)] = 1 = [G, O_3(G)]$. Fischer then proved the following.

\begin{theorem}[\cite{f}]\label{3trans}
Let $(G,D)$ be a connected irreducible $3$-transposition group.  Then up to the centre, $D$ is one of
\begin{enumerate}
\item the transposition class in $S_n$,
\item the transvection class in $O^\pm_n(2)$,
\item the transvection class in $Sp_{2n}(2)$,
\item a reflection class in $O^\pm_n(3)$,
\item the transvection class in $U_n(2)$,
\item a unique class of involutions in one of: $\Omega^+_8(2){:}S_3$, $\Omega^+_8(3){:}S_3$, $Fi_{22}$, $Fi_{23}$, or $Fi_{24}$.
\end{enumerate}
\end{theorem}

Cuypers and Hall \cite{ch} extended this by classifying the (connected) non-irreducible $3$-transposition groups and they fall into twelve cases (a convenient description can be found in \cite[Theorem 5.3]{hs}).

\subsection{Matsuo algebras}

Recall from Example \ref{Matsuo algebra} that given any $3$-transposition group $(G,D)$ and $\eta \in \FF - \{ 0,1\}$, the \emph{Matsuo algebra} $A := M_\eta(G,D)$ has basis $D$ and multiplication given by
\[
a \circ b = \begin{cases}
a & \mbox{if } b= a, \\
0 & \mbox{if } |ab| =2, \\
\frac{\eta}{2}(a+b-c) & \mbox{if }  |ab| =3, \mbox{ where } c = a^b = b^a,
\end{cases}
\]
This is an axial algebra of Jordan type $\eta$.  Moreover it also has some nice properties.  All Matsuo algebras have a Frobenius form which is given by
\[
(a,b) = \begin{cases}
1 & \mbox{if } b= a, \\
0 & \mbox{if } |ab| =2, \\
\eta/2 & \mbox{if }  |ab| =3
\end{cases}
\]
The Miyamoto group of an Matsuo algebra $M_\eta(G,D)$ is $G/Z(G)$.  In particular, we have examples of axial algebras for each group in Theorem \ref{3trans}.

Every Matsuo algebra $M_\eta(G,D)$ is an axial algebra of Jordan type $\eta$ and in fact, for $\eta \neq \frac{1}{2}$, we have the converse.

\begin{theorem}[\cite{hrs2}]\label{classification jordan type}
Every axial algebra of Jordan type $\eta \neq \frac{1}{2}$ is isomorphic to a quotient of a Matsuo algebra.
\end{theorem}

\subsection{Classification}

For $\eta = \frac{1}{2}$, Matsuo algebras are also examples of axial algebras of Jordan type $\frac{1}{2}$.  However, it is well known that a Jordan algebra satisfies the Peirce decomposition with respect to every idempotent.  This decomposition is nothing other than the statement that the idempotents are axes of Jordan type $\frac{1}{2}$.  So every Jordan algebra generated by primitive idempotents (such as say a spin factor from Example \ref{spin factor}) is an axial algebra of Jordan type $\frac{1}{2}$.

\begin{conjecture}
Every axial algebra $A$ of Jordan type $\eta$ is either isomorphic to a quotient of a Matsuo algebra, or $\eta =\frac{1}{2}$ and $A$ is a Jordan algebra.
\end{conjecture}

The conjecture is known to hold in some cases.  The $2$-generated case, that is where the algebra is generated by two primitive axes has been classified in the same paper \cite{hrs2} where the class of Jordan type algebras was introduced.  When $\eta \neq \frac{1}{2}$, we only have two Matsuo algebras.  One of these is $2\B := M_\eta(V_4, \{a,b\})$, where $V_4 = \la a,b \ra$, which is just isomorphic to $\FF \oplus \FF$ and the other is $3\C(\eta) := M_\eta(S_3, D)$, where $D$ is the class of involutions in $S_3$.  This fact almost immediately implies Theorem \ref{classification jordan type} and hence completes the classification for $\eta \neq \frac{1}{2}$.

For $\eta = \frac{1}{2}$, there is an infinite $1$-parameter family of $3$-dimensional algebras and every $2$-generated algebra is one of these, or a quotient.   The $3$-dimensional examples are, with one exception, isomorphic to the spin factor algebras $S(b)$, with $b$ a bilinear form on a $2$-dimensional vector space, as described in Example \ref{spin factor}.  The exceptional algebra doesn't have an identity and is baric; that is, it has an algebra homomorphism to $\FF$ with respect to which both generating axes map to $1$.  Almost all of these $3$-dimensional algebras are simple and so there are only quotients for some very specific parameters.  

In addition to Theorem \ref{classification jordan type}, the classification of the $2$-generated case also gives the following.

\begin{theorem}[\cite{hss2}]
Every axial algebra of Jordan type $\eta$ has a Frobenius form with respect to which every primitive axis is of length $1$.
\end{theorem}

The $3$-generated case has also been classified by Gorshkov and Staroletov \cite{gs}, namely every axial algebra of Jordan type $\frac{1}{2}$ is covered by a $9$-dimensional Jordan algebra from a $4$-parameter family.  Generically, they are isomorphic to the $3 \times 3$ matrix algebra.

The case with $4$, or more generators remains wide open and in the next Section we will discuss some recent ideas.  Note that the $4$-generated case is the first one where some Matsuo algebras are not Jordan.

\subsection{Solid subalgebras}\label{sec:solid subalgebras}

Suppose $A$ is an algebra of Jordan type $\frac{1}{2}$ and $B = \lla a,b \rra$ is a $2$-generated subalgebra.

\begin{definition}
We say that $B$ is \emph{solid} if every primitive idempotent in $B$ is a primitive axis in $A$.
\end{definition}

Every $2$-generated subalgebra of a Jordan algebra is solid.  The workshop `Algebras of Jordan type and groups'\footnote{This was supported by a Focused Research Workshop grant from the Heilbronn Institute.}, which took place in Birmingham in January 2022, was focused on the following recent result.

\begin{theorem}[Gorshkov, Shpectorov, Staroletov, 2021]
If $B$ is a spin factor and the order of $\tau_a\tau_b$ is at least $4$, then $B$ is a solid subalgebra.
\end{theorem}

This in turn implies the following in characteristic $0$.

\begin{corollary}
Every $2$-generated subalgebra $B$ of $A$ either contains $2$ primitive axes, or $3$ primitive axes, or it is solid.
\end{corollary}

It follows in the same way that Theorem \ref{classification jordan type} follows from the $2$-generated case, that if $A$ contains no solid subalgebras, then $A$ is a quotient of a Matsuo algebra.  It is tempting to conjecture that solid subalgebras identify the case where the algebra has to be Jordan.  However, the following example shows that this is not the case.

\begin{example}[Gorshkov, Staroletov, 2022]
The Matsuo algebra $M = M_{\frac{1}{2}}(3^2{:}S_3, D)$, contains exactly $3$ solid subalgebras while all other $2$-generated subalgebras contain $3$ primitive axes.  Needless to say, this means that $M$ is not a Jordan algebra\footnote{De Medts and Rehren determined in \cite{dr}, with a correction by Yabe \cite{yabe jordan}, all Matsuo algebras for $\eta = \frac{1}{2}$ that are Jordan.}.
\end{example}

So this idea seems to lead to a geometric point of view, where the primitive axes are points and we have several different types of lines: with $2$ points, $3$ points and solid lines.

\section{$2$-generated axial algebras of Monster type}
An axial algebra is of \emph{Monster type $(\al, \bt)$}, or just of Monster type, if it has the $\cM(\al, \bt)$ fusion law (see Table \ref{tab: fusion laws}).  This class is particularly interesting as it is a natural extension of the class of algebras of Jordan type.  Indeed every algebra of Jordan type $\al$, or $\bt$ is at the same time an algebra of Monster type $(\al,\bt)$.  In addition, it is also interesting as the prize example of the Griess algebra occurs here and so we are guaranteed to have algebras with many interesting groups.

We would like to classify all such axial algebras, however, as hinted at by the Griess algebra example, this is a much more difficult problem.  As with axial algebras of Jordan type, we can begin by classifying the $2$-generated algebras.

\subsection{Sakuma's Theorem}

As we have seen, the Griess algebra is a motivating example of an axial algebra and it has the fusion law $\cM(\frac{1}{4}, \frac{1}{32})$.  Its axes are in bijection with the involutions in the Monster of class $2\A$ and so they became known as $2\A$-axes.  It was Norton \cite{C85} who first studied the $2$-generated subalgebras of the Griess algebra.  He showed that the isomorphism class of the subalgebra generated by distinct $2\A$-axes $a$ and $b$ is fully determined by the conjugacy class of the element $\tau_a \tau_b$ in the Monster. There are eight such conjugacy classes of products, namely $2\A$, $2\B$, $3\A$, $3\C$, $4\A$, $4\B$, $5\A$, and $6\A$.  So, we label the $2$-generated subalgebras by these conjugacy classes.

The concept of a $2\A$-axis was generalised in the VOA setting to the concept of an Ising vector in a general OZ VOA. Miyamoto noticed in \cite{m} that the fusion property of Ising vectors gave a $\tau$ automorphism of the VOA and he proposed to study VOAs generated by Ising vectors, starting with two.  He himself finished one of the cases and the remaining cases were completed by Sakuma \cite{s}.  Amazingly, he showed that the weight $2$ component algebra of such a VOA is isomorphic to one of the above eight subalgebras of the Griess algebra.

Ivanov introduced the class of Majorana algebras \cite{i} by taking as axioms the properties he extracted from Sakuma's proof.  In particular, this is the first time that the fusion law, namely $\cM(\frac{1}{4}, \frac{1}{32})$, appeared as an axiom.  One of the early successes was the proof of Sakuma's Theorem in this setting, namely it was shown in \cite{ipss} that every $2$-generated Majorana algebra is one of the same eight Norton-Sakuma algebras. 

Axial algebras were introduced as a broad generalisation of Majorana algebras, obtained by generalising some and dropping other axioms.  Accordingly, we have the task of classifying $2$-generated algebras for broad classes of axial algebras and we call such results Sakuma-like theorems.  The first of these was for axial algebras of Monster type $(\frac{1}{4}, \frac{1}{32})$ that have a Frobenius form \cite{hrs1}.  The ultimate form of a Sakuma theorem for the class of Monster type $(\frac{1}{4}, \frac{1}{32})$ was finally obtained later as a consequence of work in \cite{generic}.

\begin{theorem}\label{sakuma thm}
Let $A$ be a $2$-generated axial algebra of Monster type $(\frac{1}{4}, \frac{1}{32})$ over a field $\FF$ of characteristic $0$.  Then, $A$ is one of the \emph{Norton-Sakuma algebras}: $2\A$, $2\B$, $3\A$, $3\C$, $4\A$, $4\B$, $5\A$, and $6\A$.
\end{theorem}

In Table \ref{tab:sakuma}, we give the structure constants and Frobenius forms for the Norton-Sakuma algebras.  For an algebra labelled $n\L$, $X = \{a_0, a_1\}$ is the set of generating axes. The Miyamoto group $\Miy(X) = \la \tau_{a_0}, \tau_{a_1} \ra$ is a dihedral group and we let $\rho := \tau_{a_0} \tau_{a_1}$.  Let $\bar{X}$ be the closure of $X$; this has size $n$.  The axes $a_i \in \bar{X}$ are labelled
\[
a_{\epsilon + 2k} := a_\epsilon^{\rho^k}
\]
where $\epsilon = 0, 1$.  It is clear from this, that every $a_i$ is conjugate to either $a_0$, or $a_1$.  In fact, if $n$ is odd, then $a_0^{\Miy(X)} = a_1^{\Miy(X)}$ and so there is a single orbit of axes under the action of the Miyamoto group.  If $n$ is even, then $a_0^{\Miy(X)}$ and $a_1^{\Miy(X)}$ are disjoint orbits and they both have size $\frac{n}{2}$.  The additional basis elements in Table \ref{tab:sakuma} which occur for some algebras are indexed by a power $\rho^k$ of $\rho$ and this indicates that the basis element is fixed by $\rho^k$.

\begin{table}[p]
\setlength{\tabcolsep}{4pt}
\renewcommand{\arraystretch}{1.5}
\centering
\footnotesize
\begin{tabular}{c|c|c}
\hline \hline
Type & Basis & Products \& form \\ \hline
$2\textrm{A}$ & \begin{tabular}[t]{c} $a_0$, $a_1$, \\ $a_\rho$ 
\end{tabular} & 
\begin{tabular}[t]{c}
$a_0 \cdot a_1 = \frac{1}{8}(a_0 + a_1 - a_\rho)$ \\
$a_0 \cdot a_\rho = \frac{1}{8}(a_0 + a_\rho - a_1)$ \\
$(a_0, a_1) = (a_0, a_\rho)= (a_1, a_\rho) = \frac{1}{8}$
\vspace{4pt}
\end{tabular}
\\
$2\textrm{B}$ & $a_0$, $a_1$ &
\begin{tabular}[t]{c}
$a_0 \cdot a_1 = 0$ \\
$(a_0, a_1) = 0$
\vspace{4pt}
\end{tabular}
\\
$3\textrm{A}$ & \begin{tabular}[t]{c} $a_{-1}$, $a_0$, \\ $a_1$, $u_\rho$ \end{tabular} &
\begin{tabular}[t]{c}
$a_0 \cdot a_1 = \frac{1}{2^5}(2a_0 + 2a_1 + a_{-1}) - \frac{3^3\cdot5}{2^{11}} u_\rho$ \\
$a_0 \cdot u_\rho = \frac{1}{3^2}(2a_0 - a_1 - a_{-1}) + \frac{5}{2^{5}} u_\rho$ \\
$u_\rho \cdot u_\rho = u_\rho$, $(a_0, a_1) = \frac{13}{2^8}$ \\
$(a_0, u_\rho) = \frac{1}{4}$, $(u_\rho, u_\rho) = \frac{2^3}{5}$ 
\vspace{4pt}
\end{tabular}
\\
$3\textrm{C}$ & \begin{tabular}[t]{c} $a_{-1}$, $a_0$, \\ $a_1$ \end{tabular} &
\begin{tabular}[t]{c}
$a_0 \cdot a_1 = \frac{1}{2^6}(a_0 + a_1 - a_{-1})$ \\
$(a_0, a_1) = \frac{1}{2^6}$
\vspace{4pt}
\end{tabular}
\\
$4\textrm{A}$ & \begin{tabular}[t]{c} $a_{-1}$, $a_0$, \\ $a_1$, $a_2$ \\ $v_\rho$ \end{tabular} &
\begin{tabular}[t]{c}
$a_0 \cdot a_1 = \frac{1}{2^6}(3a_0 + 3a_1 + a_{-1} + a_2 - 3v_\rho)$ \\
$a_0 \cdot v_\rho = \frac{1}{2^4}(5a_0 - 2a_1 - a_2 - 2a_{-1} + 3v_\rho)$ \\
$v_\rho \cdot v_\rho = v_\rho$, $a_0 \cdot a_2 = 0$ \\
$(a_0, a_1) = \frac{1}{2^5}$, $(a_0, a_2) = 0$\\
$(a_0, v_\rho) = \frac{3}{2^3}$, $(v_\rho, v_\rho) = 2$
\vspace{4pt}
\end{tabular}
\\
$4\textrm{B}$ & \begin{tabular}[t]{c} $a_{-1}$, $a_0$, \\ $a_1$, $a_2$ \\ $a_{\rho^2}$ \end{tabular} &
\begin{tabular}[t]{c}
$a_0 \cdot a_1 = \frac{1}{2^6}(a_0 + a_1 - a_{-1} - a_2 + a_{\rho^2})$ \\
$a_0 \cdot a_2 = \frac{1}{2^3}(a_0 + a_2 - a_{\rho^2})$ \\
$(a_0, a_1) = \frac{1}{2^6}$, $(a_0, a_2) = (a_0, a_{\rho^2})= \frac{1}{2^3}$
\vspace{4pt}
\end{tabular}
\\
$5\textrm{A}$ & \begin{tabular}[t]{c} $a_{-2}$, $a_{-1}$,\\ $a_0$, $a_1$,\\ $a_2$, $w_\rho$ \end{tabular} &
\begin{tabular}[t]{c}
$a_0 \cdot a_1 = \frac{1}{2^7}(3a_0 + 3a_1 - a_2 - a_{-1} - a_{-2}) + w_\rho$ \\
$a_0 \cdot a_2 = \frac{1}{2^7}(3a_0 + 3a_2 - a_1 - a_{-1} - a_{-2}) - w_\rho$ \\
$a_0 \cdot w_\rho = \frac{7}{2^{12}}(a_1 + a_{-1} - a_2 - a_{-2}) + \frac{7}{2^5}w_\rho$ \\
$w_\rho \cdot w_\rho = \frac{5^2\cdot7}{2^{19}}(a_{-2} + a_{-1} + a_0 + a_1 + a_2)$ \\
$(a_0, a_1) = \frac{3}{2^7}$, $(a_0, w_\rho) = 0$, $(w_\rho, w_\rho) = \frac{5^3\cdot7}{2^{19}}$
\vspace{4pt}
\end{tabular}
\\
$6\textrm{A}$ & \begin{tabular}[t]{c} $a_{-2}$, $a_{-1}$,\\ $a_0$, $a_1$,\\ $a_2$, $a_3$ \\ $a_{\rho^3}$, $u_{\rho^2}$ \end{tabular} &
\begin{tabular}[t]{c}
$a_0 \cdot a_1 = \frac{1}{2^6}(a_0 + a_1 - a_{-2} - a_{-1} - a_2 - a_3 + a_{\rho^3}) + \frac{3^2\cdot5}{2^{11}}u_{\rho^2}$ \\
$a_0 \cdot a_2 = \frac{1}{2^5}(2a_0 + 2a_2 + a_{-2}) - \frac{3^3\cdot5}{2^{11}}u_{\rho^2}$ \\
$a_0 \cdot u_{\rho^2} = \frac{1}{3^2}(2a_0 - a_2 - a_{-2}) + \frac{5}{2^5}u_{\rho^2}$ \\
$a_0 \cdot a_3 = \frac{1}{2^3}(a_0 + a_3 - a_{\rho^3})$, $a_{\rho^3} \cdot u_{\rho^2} = 0$\\
$(a_0, a_1) = \frac{5}{2^8}$, $(a_0, a_2) = \frac{13}{2^8}$ \\
$(a_0, a_3) = \frac{1}{2^3}$, $(a_{\rho^3}, u_{\rho^2}) = 0$
\end{tabular}\\
\hline \hline
\end{tabular}
\caption{Norton-Sakuma algebras}\label{tab:sakuma}
\end{table}

It is of course desirable to generalise Theorem \ref{sakuma thm} to positive characteristics.  Looking at Table \ref{tab:sakuma}, we see that some small characteristics $2$, $3$ and $5$ might cause problems.  In addition, we should avoid $3$, $7$ and $31$, if we want to keep the eigenvalues $1, 0, \frac{1}{4}, \frac{1}{32}$ distinct.  We expect the following is true and can possibly be obtained by generalising the arguments in \cite{generic}.

\begin{conjecture}
If $\ch(\FF)$ is sufficiently large, then every $2$-generated axial algebra of Monster type $(\frac{1}{4}, \frac{1}{32})$ is isomorphic to one of the Norton-Sakuma algebras realised over $\FF$.
\end{conjecture}


\subsection{Generalised Sakuma's Theorem}

Rehren in \cite{RT, r} was the first person to attempt to generalise the above from $\cM(\frac{1}{4}, \frac{1}{32})$ to an arbitrary Monster type $(\alpha, \beta)$.  He discovered that there are special cases $\al = 2\bt$ and $\al = 4\bt$ and proved the following:

\begin{theorem}
Suppose $\ch(\FF) \neq 2$, with $1,0, \al, \bt \in \FF$ distinct.  If $\al \neq 2\bt, 4\bt$, then every $2$-generated axial algebra $A$ of Monster type $(\al,\bt)$ has an explicit spanning set of size $8$ defined in terms of the two generating axes.  In particular, the dimension of $A$ is at most $8$.
\end{theorem}

Additionally, he generalised the Norton-Sakuma algebras to $1$- and $2$-parameter families of $2$-generated algebras of Monster type $(\al, \bt)$.  For three algebras, namely $2\B$, $2\A$ and $3\C$, generalisations exist for all values of $(\alpha, \beta)$ (excluding where $1,0,\al,\bt$ are not distinct).  The algebra $3\A$ generalises for all $(\al, \bt)$, except for $\al = \frac{1}{2}$.  However, the remaining algebras generalise to families with fusion law $\cM(\alpha, \beta)$, but only when $\al$ and $\bt$ are tied by a simple polynomial relation.

\begin{figure}[ht!]
\centering
\begin{tikzpicture}[smooth, scale=5.5,domain=-0.2:1.2, every node/.style={scale=0.85}]
\fill[red!30] (-0.2,-0.2) rectangle (1.2,1.2);
\draw[white,thick] (0.5,-0.2) -- (0.5,1.2);
\draw[red!70] (1.35,0.3) node {$3\A(\al, \bt)$};

\draw[green!40!black,thick] (0.25,-0.2) -- (0.25,1.2);
\draw[green!40!black] (0.13,0.4) node {$4\A(\frac{1}{4}, \bt)$};

\draw[blue,thick] plot (\x,{(\x)^2/2});
\draw[blue] (0.7,0.13) node {$4\B(\al, \frac{\al^2}{2})$};

\draw[green!80!black,thick,domain=-0.125:1.2] plot (\x,0.625*\x-0.125);
\draw[green!80!black] (0.8,0.52) node {$5\A(\al, \frac{5\al-1}{8})$};

\draw[purple,thick, domain=-0.2:0.4764] plot (\x, {-(\x)^2/(4*(2*\x-1))});
\draw[purple] (0.5,1.1) node[anchor=west] {$6\A(\al, \frac{-\al^2}{4(2\al-1)})$};

\draw[->,white,thick] (-0.2,0) -- (1.2,0);
\draw (1.1,0) node[anchor=north] {$\al$};
\draw[->,white,thick] (0,-0.2) -- (0,1.2);
\draw (0,1.1) node[anchor=east] {$\bt$};

\draw[white,thick] (-0.2,1) -- (1.2,1);
\draw[white,thick] (1,-0.2) -- (1,1.2);
\draw[white,thick] (-0.2,-0.2) -- (1.2,1.2);

\draw (0,1) node[anchor=north east] {$1$};
\draw (1,0) node[anchor=north east] {$1$};
\draw (0.5,0) node[anchor=north east] {$\sfrac{1}{2}$};

\draw (0.8,0.9) node {${\al=\bt}$};
\draw (0.13,0.06) node {\footnotesize ${(\sfrac{1}{4},\sfrac{1}{32})}$};
\draw (0.25,0.03125) node {$*$};
\draw (0,0) node[anchor=north east] {$0$};
\end{tikzpicture}
\caption{Rehren's generalised Norton-Sakuma algebras}
\label{fig:Rehren}
\end{figure}
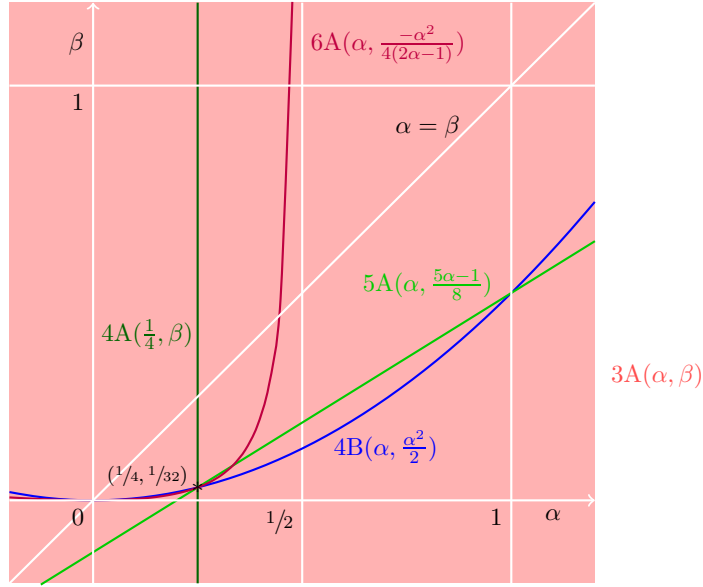

As can be seen in Figure \ref{fig:Rehren}, the point $(\frac{1}{4}, \frac{1}{32})$ is the only point in the $(\al, \bt)$-plane for which all eight examples exist.  This underscores the exceptional nature of the Griess algebra and the Monster.

\subsection{The generic case}

Rehren's work was revisited and extended by Franchi, Mainardis and Shpectorov in \cite{generic, bt 2bt, highwater}.  In \cite{generic}, they develop the general setup and derive several equations on the parameters involved.  These equations are quite complex, however it was possible to solve them in some specific situations.  For $(\al,\bt) = (\frac{1}{4}, \frac{1}{32})$, these equations yield Theorem \ref{sakuma thm}.  The other case considered is the \emph{generic case}, that is where $\al$ and $\bt$ are algebraically independent.  So in other words, we can thing of the ground field as $\QQ(\al, \bt)$, where $\al$ and $\bt$ are indeterminates.  Any algebra arising over $\QQ(\al, \bt)$ can be transferred into any characteristic and $\al$ and $\bt$ can be specialised (nearly) arbitrarily. In other words, these are the algebras that always exist.

\begin{theorem}
The generic $2$-generated algebras of Monster type are $2\B$, $3\C(\al)$, $3\C(\bt)$ and $3\A(\al, \bt)$.
\end{theorem}

As a consequence of this, the Miyamoto group of a generic algebra with an arbitrary number of generators is always a $3$-transposition group.  It is an interesting problem to determine which $3$-transposition groups actually arise here.

Additionally, Rehren's dimension bound was extended in the same paper to the exceptional cases $\al=2\bt$ and $\al=4\bt$.  In the first case, the bound is also $8$ but with a different spanning set.  In the follow up paper \cite{bt 2bt}, all the algebras in the $\al=2\bt$ case were completely classified.  The list includes the generic algebras specialised to $(2\bt,\bt)$, Joshi's algebras obtained by the double axis construction (see Section \ref{sec:double axis}) as well as a handful of exceptional algebras, some of which exist only in a specific characteristic.

In the second case, $\al=4\bt$, the situation is even more interesting.  If $(\al, \bt) \neq (2, \frac{1}{2})$, then the dimension bound of $8$ also holds, under a mild symmetry assumption (see Section \ref{sec:symmetric 2-gen}).  However, for $(\al, \bt) = (2, \frac{1}{2})$, there is an infinite dimensional example, called the Highwater algebra.

\subsection{The Highwater algebra}

The Highwater algebra $\cH$ was discovered independently by Franchi, Mainardis and Shpectorov \cite{highwater}, and by Yabe \cite{yabe}, and it is an infinite dimensional algebra of Monster type $(2, \frac{1}{2})$.  It has infinitely many axes $a_i$, one for each integer $i \in \ZZ$, and these can be completed to a basis of $\cH$ by adding $s_j$, for $j \in \NN$, which correspond to the `distance' between the axes.

\begin{definition}
Let $\ch(\FF) \neq 2$.  The \emph{Highwater algebra} is the algebra $\cH$ on $\bigoplus_{i \in \ZZ} \FF a_i \oplus \bigoplus_{j \in \NN} \FF s_j$ with multiplication given by
\begin{align*}
a_i a_j &:= \tfrac{1}{2}(a_i+a_j) +s_{|i-j|} \\
a_i  s_{j} &:= -\tfrac{3}{4} a_i + \tfrac{3}{8}( a_{i-j}+ a_{i+j}) +\tfrac{3}{2} s_{j} \\
s_j  s_ k &:= \tfrac{3}{4}( s_{j}+ s_{k}) - \tfrac{3}{8}(s_{|j-k|} + s_{j+k})
\end{align*}
for all $i \in \ZZ$ and $j,k \in \NN$, where we set $s_0 = 0$.
\end{definition}

In characteristic $3$, it requires infinitely many generators; otherwise it is $2$-generated. The algebra is also baric, that is the map $\lm \colon A \to \FF$ defined by $\lm(a_i) = 1$ and $\lm(s_j) = 0$ is an algebra homomorphism.  This immediately shows that $\cH$ has a Frobenius form, given by $(x,y) = \lm(x)\lm(y)$, for all $x,y \in \cH$.

For every axis $a_i \in \cH$, the Miyamoto involution $\tau_i := \tau_{a_i}$ is given by $a_j \mapsto a_{2i-j}$ which is a reflection on the indices of the axis.  It is then easy to see that $\Miy(\cH) \cong D_\infty$.  There is also an additional symmetry given by swapping $a_0$ and $a_1$, which we write as $\tau_{\frac{1}{2}}$ and is a reflection in $\frac{1}{2}$.  It turns out that $\Aut(\cH) = \la \Miy(\cH), \tau_{\frac{1}{2}} \ra \cong D_\infty$.

In characteristic 5 only, the Highwater algebra admits a cover $\hatH$ found by Franchi and Mainardis \cite{highwater5}.  This is also a $2$-generated axial algebra of Monster type $(2, \frac{1}{2})$.  It is slightly more complicated than $\cH$, with extra basis elements $p_{\bar{r},j}$, where $\bar{r} \in \{\bar{1}, \bar{2}\} \in \ZZ_3$ and $j \in 3\NN$ (see \cite[Definition 3.2]{HWquo}).  These generate an ideal $J = \la p_{\bar{r},j} : \bar{r} \in \{\bar{1}, \bar{2}\}, j \in 3\NN \ra$ and $\hatH/J \cong \cH$.


\subsection{Symmetric $2$-generated algebras of Monster type}\label{sec:symmetric 2-gen}

One can note that all of the (generalised) Norton-Sakuma algebras have the property that there is an involutory automorphism which switches the two generating axes.  We call such algebras \emph{symmetric}.  In a major breakthrough, Yabe \cite{yabe} produced an almost complete classification of the symmetric $2$-generated algebras of general Monster type $(\alpha, \beta)$.  As part of this, he discovered some further infinite families of examples in addition to Rehren's generalised Norton-Sakuma algebras and Joshi's examples coming from the double axis construction (which we will see in Section \ref{sec:double axis}) and he also independently discovered the Highwater algebra.  The missing characteristic $5$ case was completed by Franchi and Mainardis \cite{highwater5}, who showed that the only additional examples are quotients of the cover $\hatH$ of the Highwater algebra.  Finally, the quotients of the Highwater algebra and its cover were classified by Franchi, Mainardis and M\textsuperscript{c}Inroy in \cite{HWquo}.

\begin{theorem}[\cite{yabe, highwater5, HWquo}]\label{2gen sym}
A symmetric $2$-generated axial algebra of Monster type $(\alpha, \beta)$ is a quotient of one of the following:
\begin{enumerate}
\item an axial algebra of Jordan type $\alpha$, or $\beta$;
\item in one of the following families:
\begin{description}
\item[\textup{\normalfont(a)}] $3\A(\al,\bt)$, $4\A(\frac{1}{4}, \bt)$, $4\B(\al, \frac{\al^2}{2})$, $4\J(2\bt, \bt)$, $4\Y(\frac{1}{2}, \bt)$, $4\Y(\al, \frac{1-\al^2}{2})$, $5\A(\al, \frac{5\al-1}{8})$, \linebreak $6\A(\al, \frac{-\al^2}{4(2\al-1)})$, $6\J(2\bt, \bt)$ and $6\Y(\frac{1}{2}, 2)$;
\item[\textup{\normalfont(b)}] $\IY_3(\al, \frac{1}{2}, \mu)$ and $\IY_5(\al, \frac{1}{2})$;
\end{description}
\item[$3$] the Highwater algebra $\cH$, or its characteristic $5$ cover $\hatH$.
\end{enumerate}
\end{theorem}

Note that our cases and notation are different to Yabe's\footnote{Yabe's names for these algebras are: (a) $\textrm{III}(\al, \bt,0)$, $\textrm{IV}_1(\frac{1}{4}, \bt)$, $\textrm{IV}_2(\al, \frac{\al^2}{2})$, $\textrm{IV}_1(\al, \frac{\al}{2})$, $\textrm{IV}_2(\frac{1}{2}, \bt)$, $\textrm{IV}_2(\al, \frac{1-\al^2}{2})$, $\textrm{V}_1(\al, \frac{5\al-1}{8})$, $\textrm{VI}_2(\al, \frac{-\al^2}{4(2\al-1)})$, $\textrm{VI}_1(\al, \frac{\al}{2})$ and $\textrm{IV}_3(\frac{1}{2}, 2)$ and (b) $\textrm{III}(\al, \frac{1}{2}, \dl)$, where $\dl = -2\mu-1$, and $\textrm{V}_2(\al, \frac{1}{2})$.  Here the Roman numeral indicates axial dimension.}.  Similarly to the (generalised) Norton-Sakuma algebras, we label an algebra $n\L$ to emphasis the size $n$ of the closed set of axes.  The two algebras which are labelled $\IY_d$ are of axial dimension $d$ and generically have infinitely many axes, however they can have finitely many axes for some characteristics and values of their parameters.  We use the same labels $\L$ as for the (generalised) Norton-Sakuma algebras, adding $\J$ for Joshi's examples and $\Y$ for the algebras introduced by Yabe.  The parameters after this are $(\al,\bt)$ (and in the case of $\IY_3(\al, \frac{1}{2}, \mu)$ there is an additional parameter $\mu$).

\begin{table}[htb!]
\setlength{\tabcolsep}{4pt}
\renewcommand{\arraystretch}{1.4}
\centering
\footnotesize
\begin{tabular}{c|c|c|c|c|c}
\hline \hline
Name & $\{1,0,\alpha,\beta\}$ coinciding & Additional & Dimension & Quotients & Dimension\\
\hline
$3\A(\al,\bt)$ & $\alpha, \beta \neq 0,1$, $\alpha \neq \beta$ & $\alpha \neq \frac{1}{2}$ & 4 & $3\A(\al,\frac{1-3\al^2}{3\al-1})^\times$ & 3 \\
$4\A(\frac{1}{4}, \bt)$ & $\ch(\FF) \neq 3$, $\alpha \neq 0,1,\frac{1}{4}$ & & 5 & $4\A(\frac{1}{4},\frac{1}{2})^\times$ & 4 \\
$4\B(\al, \frac{\al^2}{2})$ & $\al \neq \{ 0,1,2, \pm \sqrt{2}\}$ & & 5 & $4\B(-1, \frac{1}{2})^\times$ & 4 \\
$4\J(2\bt, \bt)$ & $\bt \neq \{ 0,1,\frac{1}{2} \}$ & & 5 & $4\J(-\frac{1}{2}, -\frac{1}{4})^\times$ & 4 \\
$4\Y(\frac{1}{2}, \bt)$ & $\bt \neq \{ 0,1,\frac{1}{2} \}$ & & 5 & & \\
$4\Y(\al, \frac{1-\al^2}{2})$ &
\begin{tabular}[t]{c}
$\al \neq \{ 0,\pm 1, \pm \sqrt{-1},$ \\
\quad $ -1 \pm \sqrt{2}\}$
\end{tabular}& & 5 & & \\
$5\A(\al, \frac{5\al-1}{8})$ & $\al \neq \{ 0,1, -\frac{1}{3}, \frac{1}{5}, \frac{9}{5} \}$ & & 6 & & \\
$6\A(\al, \frac{-\al^2}{4(2\al-1)})$ &
\begin{tabular}[t]{c}
$\al \neq \{ 0, 1, \frac{4}{9}, \qquad$ \\
\hspace{30pt}$ -4 \pm 2\sqrt{5}\}$
\end{tabular}
& $\al \neq \frac{1}{2}$ & 8 & 
\begin{tabular}[t]{c}
$6\A(\frac{2}{3}, -\frac{1}{3})^\times$ \\
$6\A(\frac{1 \pm \sqrt{97}}{24}, \frac{53 \pm5\sqrt{97}}{192})^\times$
\end{tabular} &  \begin{tabular}[t]{c} 7 \\ 7 \end{tabular}\\
$6\J(2\bt, \bt)$ & $\bt \neq \{ 0,1,\frac{1}{2} \}$ & & 8 & $6\J(-\frac{2}{7}, -\frac{1}{7})^\times$ & 7\\
$6\Y(\frac{1}{2}, 2)$ & $\ch(\FF) \neq 3$ & & 5 & $6\Y(\frac{1}{2}, 2)^\times$ & 4\\
$\IY_3(\al, \frac{1}{2}, \mu)$ & $\alpha \neq 0,1,\frac{1}{2}$ & & 4 & 
\begin{tabular}[t]{c}
$\IY_3(-1, \frac{1}{2}, \mu)^\times$ \\
$\IY_3(\al, \frac{1}{2}, 1)^\times$
\end{tabular} & \begin{tabular}[t]{c} 3 \\ 3 \end{tabular}
\\
$\IY_5(\al, \frac{1}{2})$ & $\alpha \neq 0,1,\frac{1}{2}$ & & 6 & $\IY_5(\al, \frac{1}{2})^\times$ & 5 \\
\hline \hline
\end{tabular}
\caption{The symmetric $2$-generated axial algebras of Monster type}\label{tab:2-gen Monster}
\end{table}

We list some properties of the algebras in Table \ref{tab:2-gen Monster}, but we do not propose to give the structure constants for all them here.  See \cite{yabe} for the original definition of these algebras and \cite{forbidden} for new bases for some of them.  However, we will comment on the permitted values $(\alpha, \beta)$ for these algebras as there are some errors in both Rehren's and Yabe's tables.  It is clear from our definition of axial algebras, that we must have that $\{1,0,\al,\bt\}$ are distinct, otherwise some eigenspaces would merge (however, see decomposition algebras \cite{dpsv} for a way round this).  For all the algebras except $3\A(\al, \bt)$ and $6\A(\al, \frac{-\al^2}{4(2\al-1)})$, this is the only restriction.  For these two, we additionally only need to exclude $\al = \frac{1}{2}$ (Rehren and Yabe excluded other values here unnecessarily).  This information can be readily checked on the computer in many cases.  M\textsuperscript{c}Inroy has written a Magma package \cite{2-gen magma} with these algebras and quotients of the Highwater algebra.

\begin{question}
The value $\al=\frac{1}{2}$ is excluded for $3\A(\al, \bt)$ as some of the structure constants are not well-defined.  Is there another description of this family which makes sense for $\al=\frac{1}{2}$?  Note that for $6\A(\al, \frac{-\al^2}{4(2\al-1)})$, $\beta$ tends to infinity as $\al$ tends to $\frac{1}{2}$, so this cannot be helped.
\end{question}

It is also evident from Table \ref{tab:2-gen Monster}, that any symmetric $2$-generated algebra with $\beta \neq \frac{1}{2}$ has at most $6$ axes.  It turns out that when $\bt = \frac{1}{2}$, the algebras $\IY_3(\al, \frac{1}{2}, \mu)$, $\IY_5(\al, \frac{1}{2})$, $\cH$ and $\hatH$ have generically infinitely many axes and they (or their quotients) can have any finite number of axes.

\begin{problem}
Explain the difference in the number of axes when $\bt = \frac{1}{2}$ and $\beta \neq \frac{1}{2}$.
\end{problem}

Possibly linked to this, Krasnov and Tkachev \cite{kv} have found exceptional behaviour for the eigenvalue $\frac{1}{2}$.  Roughly speaking if $\frac{1}{2}$ appears as an eigenvalue of any semisimple idempotent, then the total number of idempotents is at most $2^{\dim(A)}$.

We note that every symmetric $2$-generated algebra of Monster type in Theorem \ref{2gen sym} has a Frobenius form.

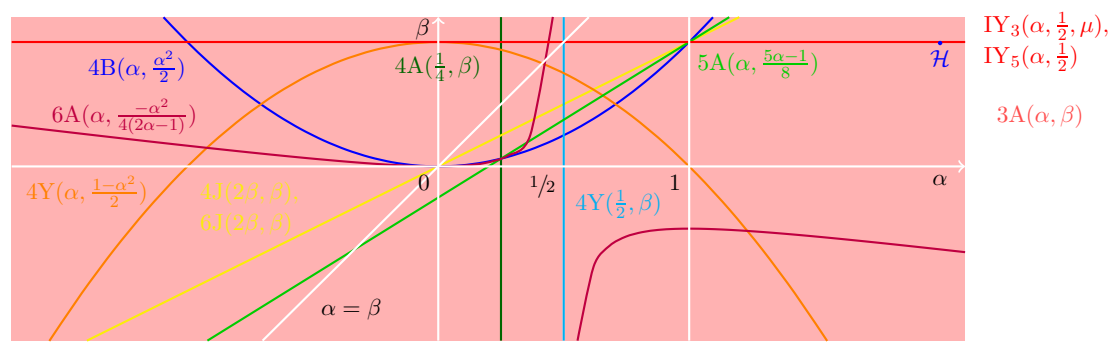
\begin{figure}[ht!]
\centering
\begin{tikzpicture}[smooth, scale=3.3, yscale=1, domain=-1.7:2.1, every node/.style={scale=0.8}]
\fill[red!30] (-1.7,-0.7) rectangle (2.1,0.6);
\draw[white,thick] (0.5,-0.7) -- (0.5,0.6);
\draw[red!70] (2.4,0.2) node {$3\A(\al, \bt)$};

\draw[green!40!black,thick] (0.25,-0.7) -- (0.25,0.6);

\draw[blue,thick, domain=-1.095:1.095] plot (\x,{(\x)^2/2});
\draw[blue] (-1.2,0.4) node {$4\B(\al, \frac{\al^2}{2})$};

\draw[yellow!98!black, thick, domain=-1.4:1.2] plot (\x,{\x/2});
\draw[yellow!98!black] (-0.75,-0.17) node {$\begin{array}{l} 4\J(2\bt, \bt), \vspace{1pt}\\ 6\J(2\bt, \bt) \end{array}$};

\draw[cyan,thick] (0.5,-0.7) -- (0.5,0.6);
\draw[cyan] (0.72,-0.15) node {$4\Y(\frac{1}{2}, \bt)$};

\draw[orange,thick, domain=-1.55:1.55] plot (\x,{(1-(\x)^2)/2});
\draw[orange] (-1.4,-0.1) node {$4\Y(\al, \frac{1-\al^2}{2})$};

\draw[green!80!black,thick,domain=-0.92:1.16] plot (\x,0.625*\x-0.125);
\draw[green!80!black] (1,0.5) node[anchor=north west] {$5\A(\al, \frac{5\al-1}{8})$};

\draw[purple,thick, domain=-1.7:0.45662] plot (\x, {-(\x)^2/(4*(2*\x-1))});
\draw[purple,thick, domain=0.555:2.1] plot (\x, {-(\x)^2/(4*(2*\x-1))});
\draw[purple] (-1.25,0.2) node {$6\A(\al, \frac{-\al^2}{4(2\al-1)})$};


\draw[red,thick] (-1.7,0.5) -- (2.1,0.5);
 \draw[red] (2.1,0.5) node[anchor=west] {$\begin{array}{l} \IY_3(\al, \frac{1}{2}, \mu), \vspace{1pt} \\ \IY_5(\al, \frac{1}{2}) \end{array}$};

\draw[blue] (2,0.505) node[circle, fill,scale=.2,anchor=north] {};
\draw[blue] (2,0.5) node[anchor=north] {$\cH$};

\draw[->,white,thick] (-1.7,0) -- (2.1,0);
\draw (2,0) node[anchor=north] {$\al$};
\draw[->,white,thick] (0,-0.7) -- (0,0.6);
\draw (0,0.55) node[anchor=east] {$\bt$};

\draw[green!40!black] (0,0.4) node {$4\A(\frac{1}{4}, \bt)$};

\draw[white,thick] (1,-0.7) -- (1,0.6);
\draw[white,thick] (-0.7,-0.7) -- (0.6,0.6);

\draw (1,0) node[anchor=north east] {$1$};
\draw (0.5,0) node[anchor=north east] {$\sfrac{1}{2}$};

\draw (-0.5,-0.5) node[anchor =north west] {${\al=\bt}$};
\draw (0,0) node[anchor=north east] {$0$};
\end{tikzpicture}
\caption{Symmetric $2$-generated axial algebras of Monster type $(\al,\bt)$}
\label{fig:Yabe}
\end{figure}

Finally, in Figure \ref{fig:Yabe}, we give a picture similar to Figure \ref{fig:Rehren} showing the intersection of the permitted values of $(\al, \bt)$ for the symmetric $2$-generated algebras of Monster type.  From this picture, we can now see some additional interesting points where several algebras exist at the same time.

\begin{problem}
Investigate the algebras and associated Miyamoto groups occurring at the intersection points in Figure \ref{fig:Yabe}.
\end{problem}

\subsection{Quotients of the Highwater algebra and its cover}

In \cite{HWquo}, Franchi, Mainardis and M\textsuperscript{c}Inroy classify the ideals of the Highwater algebra and its cover.  In fact, they introduce a cover $\hatH$ in all characteristics, which is necessarily not of Monster type except in characteristic $5$.  They then use this to produce a unified classification for the ideals, which can then be specialised to the Highwater algebra, or its cover in characteristic $5$.

\begin{theorem}
Every ideal of the Highwater algebra $\cH$ and its cover $\hatH$ is invariant under the full automorphism group.  In particular, every quotient is a symmetric $2$-generated algebra of Monster type $(2, \frac{1}{2})$.
\end{theorem}

For $\cH$, they show that every ideal has finite codimension -- so there are no infinite dimensional quotients of $\cH$.  In particular, we say that an ideal $I$ has \emph{axial codimension} $D$ if the quotient $\cH/I$ has axial dimension $D$.  That is, if the images of the axes in quotient span a $D$-dimensional subspace.  The results can now be stated very succinctly with the use of one definition.

\begin{definition}
A tuple $(\al_0, \dots, \al_D) \in \FF^{D+1}$ is of \emph{ideal type} if
\begin{enumerate}
\item $\al_0 \neq 0 \neq \al_D$
\item $\sum_{i = 0}^D \al_i = 0$
\item There exists $\epsilon = \pm 1$, such that for all $i = 0, \dots, D$, $\al_i = \al_{D-i}$.
\end{enumerate}
\end{definition}

\begin{theorem}[\cite{HWquo}]\label{highwater quo}
For each $D \in \NN$, there is a bijection between the ideal type tuples, up to scalars, and the minimal ideals of $\cH$ of axial codimension $D$ given by
\[
(\al_0, \dots, \al_D) \mapsto \left( \sum_{i=0}^D \al_i a_i \right) \unlhd \cH
\]
\end{theorem}

In particular, it is clear that $(a_0-a_D) \unlhd \cH$, for all $D \in \NN$, and hence $\cH$ has a quotient with $D$ axes for every $D \in \NN$.  In addition, the above theorem also leads to an explicit basis for a general ideal and hence for the quotient.  

The ideals for $\hatH$ are only slightly more complicated.  Recall that for $\hatH$, there is a distinguished ideal $J$ such that $\hatH/J \cong \cH$.  Every ideal $I \unlhd \hatH$ has either finite codimension in $\hatH$ and is described by the above theorem, or it is contained in the ideal $J$ and has finite codimension in $J$.  The paper \cite{HWquo} contains a theorem on the ideals contained in $J$ and explicit bases, that is similar to Theorem \ref{highwater quo}.


\subsection{Non-symmetric $2$-generated axial algebras}

The case of non-symmetric $2$-generated algebras of Monster type $(\al,\bt)$ remains open.  However, some examples are known which we will briefly discuss.

\begin{example}\label{ex:Q2}
One of the series of the algebras built by Joshi using the double axis construction (see Section \ref{sec:double axis}) is called $Q_2(2\bt,\bt)$ and this is a $1$-parameter family of $2$-generated algebras of Monster type $(2\bt, \bt)$.  It contains four axes in total, two of which are single axes and the other two are double axes.  The algebra is non-symmetric as single and double axes have different lengths with respect to the unique Frobenius form.
\end{example}

The other known non-symmetric examples come from a generalised Norton-Sakuma algebra, which was noticed by M\textsuperscript{c}Inroy.

\begin{example}
For $\al = -1$, the algebra $A = 4\B(-1, \frac{1}{2})$ has a $2$-dimensional radical $R \unlhd A$.  Furthermore, every $1$-dimensional subspace of $R$ is also an ideal.  It can be seen that the automorphism $\sigma$ switching the two generating axes $a_0$ and $a_1$ of $A$ acts on $R$ with eigenvalues $1$ and $-1$.  Hence every $1$-dimensional space in $R$ that is not one of the two eigenspaces is not invariant under $\sigma$ and hence it gives a quotient of $A$ that is not a symmetric algebra.
\end{example}

Let's call a $2$-generated algebra $A = \lla a, b \rra$ \emph{skew} if $a^{\Miy(A)}$ and $b^{\Miy(A)}$ have different lengths.  The two examples above are not skew.  It is clear that a skew algebra cannot be symmetric.  This idea ties into the concept of an axet which we will discuss in Section \ref{sec:axets} and we will see there an example of a skew $2$-generated algebra of Monster type.

\begin{problem}
Classify the non-symmetric $2$-generated axial algebras of Monster type.
\end{problem}

\section{Larger algebras}
Beyond the $2$-generated case, very little is known about axial algebras of Monster type $(\al,\bt)$.  In particular, we are far from knowing all the examples and hence not in a position to formulate a conjecture for the complete list.  However, we do have some recent constructions of examples, both theoretical and computational, that demonstrate that there is a very rich variety of such algebras.

We note that all of the currently known examples of algebras of Monster type admit a Frobenius form.  Hence we can at least formulate the following conjecture.

\begin{conjecture}
Every axial algebra of Monster type $(\al, \bt)$ admits a Frobenius form.
\end{conjecture}

Let us now turn to some families of examples.


\subsection{Split spin factor algebras}

An interesting new family of algebras of Monster type arose when trying to understand the symmetric $2$-generated algebras discovered by Yabe. In particular, M\textsuperscript{c}Inroy and Shpectorov found in his algebra $\textrm{III}(\al,\frac{1}{2},\dl)$ (denoted in Theorem \ref{2gen sym} as $\IY_3(\al, \frac{1}{2}, \mu)$) a different, more natural basis, and the way multiplication works in this basis readily generalises to any number of generators. 

\begin{definition}[\cite{splitspin}]
Let $\FF$ be a field of characteristic other than $2$ and let $E$ be a quadratic space over $\FF$, i.e., a vector space of $\FF$ endowed with a symmetric bilinear form $b \colon E\times E\to\FF$. Choose $\al\in\FF\setminus\{0,\frac{1}{2},1\}$. The \emph{split spin factor algebra} is the commutative algebra $S(b, \al)=\FF z_1\oplus\FF z_2\oplus E$ with multiplication given by
\[
\begin{gathered}
z_1^2=z_1, \quad z_2^2=z_2, \quad z_1z_2=0,\\
z_1e=\al e, \quad z_2 e=(1-\al)e,\\
ef=-b(e,f)z,
\end{gathered}
\]
where $z=(\al-2)\al z_1+(\al-1)(\al+1)z_2$ and $e,f\in E$.
\end{definition} 

It is clear that $\1 = z_1+z_2$ is the identity of the algebra. The above construction is very reminiscent of the construction of the spin factor algebra in Example \ref{spin factor}. There one extends a quadratic space by a $1$-dimensional algebra $1A\cong\FF$ containing the identity; here we instead extend $E$ by a $2$-dimensional algebra $2B\cong\FF\oplus\FF$ in such a way that the generating idempotents split the identity and act on $E$ asymmetrically by different scalars.

The nonzero idempotents in $A = S(b,\al)$ can be listed explicitly. In addition to the obvious idempotents $z_1$, $z_2$, and $\1$, there are only two families of idempotents: 
\begin{enumerate}
\item[(a)] $\frac{1}{2}(e+\al z_1+(\al+1)z_2)$,
\item[(b)] $\frac{1}{2}(e+(2-\al)z_1+(1-\al)z_2)$,
\end{enumerate}
where $e\in E$ with $b(e,e)=1$.

It turns out that $z_1$ and $z_2$ are primitive idempotents of Jordan type $\al$ and $1-\al$ respectively and, similarly, the idempotents from (a) (respectively, (b)) are primitive axes of Monster type $(\al,\frac{1}{2})$ (respectively, $(1-\al,\frac{1}{2})$). The following result now follows:

\begin{theorem}
The split spin factor algebra $S(b,\al)$ is an algebra of Monster type $(\al,\frac{1}{2})$ \textup{(}and $(1-\al,\frac{1}{2})$\textup{)} if and only if $\al\neq -1$ \textup{(}respectively, $2$\textup{)} and $E$ is spanned by vectors $e$ with $b(e,e)=1$.
\end{theorem}

For the first type we select $z_1$ and family (a) as the set of generating axes and, symmetrically, for the second type we select $z_2$ and family (b).

When $\al=2$ (respectively, $-1$), $A$ is baric with the non-generator idempotents being contained in the radical. Otherwise, $A$ is simple as long as $b$ is non-degenerate. 

If $\dim E\geq 2$, the automorphism group of $A$ coincides with the orthogonal group $O(E,b)$ of the quadratic space $E$. When $\dim E=1$, the eigenvalue $\frac{1}{2}$ disappears and so $A$ is $3$-dimensional algebra of Jordan type $\al$, assuming that $E$ is panned by a vector $e$ with $b(e,e)=1$. Namely, $A$ is isomorphic to $3\C(\al)$ in this case.

Finally to relate this back to Yabe's algebras, let $E$ be a $2$-dimensional vector space with basis $\{e, f\}$ where $b(e,e) = 1 = b(f,f)$.  Note that the bilinear form is then entirely determined by $\mu := b(e,f)$.

\begin{theorem}
If $\al \neq -1$ and $\mu \neq 1$, then $S(b, \al) \cong \IY_3(\al, \frac{1}{2}, \mu)$.
\end{theorem}


\subsection{Double axis construction}\label{sec:double axis}

The double axis construction was introduced by Joshi in \cite{j, JoshiPhD} and appeared in print in a much developed form in \cite{doubleMatsuo}.

Let $M = M_\eta(G,D)$ be a Matsuo algebra, where $\eta \neq \frac{1}{2}$.  We say that the axes $a, b \in D$ are \emph{orthogonal} if $ab=0$, or equivalently $a$ and $b$ commute as involutions in $G$.  If $a$ and $b$ are orthogonal, then $x := a+b$ is an idempotent in $M$ and we will call such idempotents \emph{double axes}.

\begin{theorem}[\cite{j}]
Every double axis satisfies the fusion law $\cM(2\eta, \eta)$.
\end{theorem}

The double axis $x$ is not primitive as both $a$ and $b$ are in the $1$-eigenspace of $x$.  However, double axes may be primitive in some proper subalgebra.  Joshi in \cite{j} explored this idea and constructed all primitive algebras generated by $2$ axes, where these may be single (note that every $d \in D$ trivially satisfies $\cM(2\eta, \eta)$ and so we call these \emph{single axes}), or double axes and at least one of these is a double axis.  Joshi found three new infinite families of examples.

\begin{theorem}[\cite{j}]\label{Joshi 2-gen}
A primitive $2$-generated subalgebra of $M$, where at least one of the generating axes is a double axis, is one of $Q_2(2\eta, \eta)$, $4\J(2\eta, \eta)$, or $6\J(2\eta, \eta)$.
\end{theorem}

Note that the first of these is a non-symmetric example mentioned in Example \ref{ex:Q2} and the other two are symmetric and can be seen in the statement of Theorem \ref{2gen sym}.

It is easy to see that the two generators of the subalgebra in Theorem \ref{Joshi 2-gen} involve no more than four elements of $D$.  This means that the entire calculation happens within a small ambient group of $3$-transpositions.  The cases are organised in terms of the diagram on the set of generators of this ambient group.  It was observed that the cases leading to the examples in Theorem \ref{Joshi 2-gen} are those where the diagram admits a automorphism of order $2$ that extends to an automorphism of the ambient group.  We call this a \emph{flip}.  This leads to the following general construction.

Let $\sigma \in \Aut(G)$ of order $2$ which preserves $D$.  Then the group $S = \la \sigma \ra$ acts on $D$ and the orbits are classified into three groups: orbits of length $1$, orbits of length $2$ which are orthogonal and orbits of length $2$ which are non-orthogonal.  Taking orbit sums, we obtain single axes and double axes from the first two, and vectors called \emph{extras} from the last.  All these vectors lie in the fixed subalgebra $M^\sigma$ and in fact they form a basis.  The subalgebra of $M^\sigma$ generated by all single and double axes above is called the \emph{flip subalgebra}.

\begin{theorem}[\cite{doubleMatsuo}]
Every flip subalgebra is a \textup{(}primitive\textup{)} axial algebra of Monster type $(2\eta, \eta)$.
\end{theorem}

As can be seen from Theorem \ref{3trans}, there are many classes of groups of $3$-transpositions and they admit many different flips.  Hence flip subalgebras give a very large family of examples.

Joshi determined all flips and flip subalgebras for $G = S_n$ and $G = Sp_{2n}(2)$ \cite{j, JoshiPhD}.  Alsaeedi \cite{aPhD, apaper} completed the $2^{n-1}{:}S_n$ case.  Shi \cite{shi masters} analysed the case of $O_{2n}^\pm(2)$ and finally, in as yet unpublished work, Hoffman, Rodrigues and Shpectorov finished the $U_n(2)$ case.

\begin{question}
Are there any examples of axial algebras of Monster type $(2\eta, \eta)$ which cannot be constructed as a flip subalgebra?
\end{question}

Such an example would require at least three generators.


\subsection{Axets and shapes}\label{sec:axets}

While computing axial algebras for concrete groups, one has to deal with a closed set of generating axes together with the action of the Miyamoto group on it and the $\tau$-map.  Different configurations of $2$-generated subalgebras correspond to different \emph{shapes}, introduced by Ivanov in the context of Majorana algebras.  In essence, shapes are the cases which one needs to consider.  However, this leads to a slight technical problem: we want to talk about shapes which collapse, i.e.\ there is no algebra with that shape, so shapes should not be defined in terms of an ambient algebra.

\subsubsection{Axets}

The following concept, introduced by M\textsuperscript{c}Inroy and Shpectorov in \cite{forbidden}, frees the set of axes $X$ from the algebra.  Recall that the $\tau$-map relates the Miyamoto automorphisms to each axis.

\begin{definition}
Let $S$ be an abelian group.  An \emph{$S$-axet} $(G, X, \tau)$ is a $G$-set $X$ together with a map $\tau \colon X \times S \to G$ (written $\tau_x(s) = \tau(x,s)$), such that
\begin{enumerate}
\item $\tau_x(s) \in G_x$
\item $\tau_x(ss') = \tau_x(s) \tau_x(s')$
\item $\tau_{xg}(s) = \tau_x(s)^g$
\end{enumerate}
for all $s,s' \in S$, $x \in X$ and $g \in G$.
\end{definition}

We call elements of $X$ \emph{axes}.  Clearly, this definition mimics how, in an axial algebra, the Miyamoto group $G = \Miy(A)$ acts on the closed set of axes $X$, linked via the $\tau$-map.  In particular, every axial algebra $A$ contains an axet $X(A)$.

Note that, in this survey, we are mostly interested in $C_2$-graded algebras, so here we will take $S = C_2$ and talk about $C_2$-axets.

Even though an axet is just a combinatorial and group-theoretic object, it still carries much information about the axes in an axial algebra.  For example, the Miyamoto group of an axet $(G,X, \tau)$ is defined in the obvious way as $\Miy(X) = \la \mathrm{Im}(\tau) \ra \unlhd G$.  We can also talk about subaxets and, in particular, the subaxet $\la Y \ra$ generated by a subset $Y \subseteq X$.

\begin{example}
Let $n \in \NN \cup \{ \infty\}$ and consider the natural action of $G = D_{2n}$ on the regular $n$-gon (or $\ZZ$ if $n = \infty$).  Let $X$ be the set of vertices and $\tau_x$ be the reflection in $x$, for each $x \in X$.  We call $X(n) := (G, X, \tau)$ the \emph{$n$-gonal $C_2$-axet}.
\end{example}

It is not too difficult to see that if $n$ is odd, then there is one orbit of axes under $\Miy(X)$ and if $n$ is even, or infinite, then there are two orbits of equal size.  Moreover, if $n$ is even, then opposite vertices in the $n$-gon give the same reflection and hence the same Miyamoto involution.  This observation forms the basis of our next example.

\begin{example}
Consider the axet $X = X(4k)$.  The Miyamoto group $\Miy(X) \cong D_{4k}$ has two orbits $X_1$ and $X_2$ on $X$, each of length $2k$.  In particular, each vertex in the $4k$-gon is in the same orbit as its opposite vertex.  We now `fold' one of these orbits, $X_1$ say, by identifying opposite vertices.  Let $Y$ be the set corresponding to pairs of a vertex and its opposite vertex.  Finally, let $X' = Y \cup X_2$ and let $\tau' \colon X' \to D_{4k}$ be the obvious adjustment of the map $\tau$.  Then $X'(3k) := (D_{4k}, X', \tau')$ is a $C_2$-axet of size $3k$ with two orbits of axes, one of size $k$ and the other of size $2k$.  
\end{example}

\begin{theorem}[\cite{forbidden}]
Let $X = \la a,b\ra$ be a $2$-generated $C_2$-axet with $n$ axes, where $n \in \NN \cup \{ \infty\}$.  Then $(\Miy(X), X, \tau)$ is isomorphic \textup{(}up to the kernel of the action\textup{)} to
\begin{enumerate}
\item $X(n)$, or
\item $X'(3k)$, where $n=3k$ and $k \in \NN$.
\end{enumerate}
\end{theorem}

Note that we haven't used the fusion law, or any algebraic information here, and yet we can classify the possibilities for the axets of any $2$-generated axial algebra with an arbitrary $C_2$-graded fusion law.  Clearly, the second axet $X'(3k)$ cannot arise from a symmetric $2$-generated algebra.  In fact, all the examples of a $2$-generated algebra we have seen, have the $X(n)$ axet.  One can ask whether the \emph{skew} axet $X'(3k)$ can arise from a $2$-generated axial algebra.

The following observation by Turner gives a positive answer and hints of a possible class of very interesting algebras.

\begin{example}
Let $A = 3\C(\eta)$ with axes $a,b,c$.  If $\eta \neq -1$, then it is not difficult to see that $A$ has an identity $ \1 := \frac{1}{\eta+1}(a+b+c)$.  Now, $a' := \1-a$ is an idempotent and it turns out it is an axis of Jordan type $1-\eta$.  Consider $A$ as an axial algebra of Monster type $(\al, 1-\al)$ with the two generators $a'$ and $b$. Then, $\tau_{a'}$ switches $b$ and $c$ while $\tau_b$ and $\tau_c$ are trivial. Hence we see that $A = \lla a', b \rra$ is an axial algebra of Monster type $(\al, 1-\al)$ with a skew axet $X'(1+2)$.
\end{example}

However, this example is somewhat degenerate and no further examples are known.

\begin{problem}
Are there any other $2$-generated axial algebras of Monster type $(\al, \bt)$ with a skew axet?  If so, classify them.
\end{problem}

M\textsuperscript{c}Inroy and Shpectorov \cite{binary} also classified the \emph{binary} axets, which are the the axets where all orbits of the Miyamoto group have size $2$.  This classification is in terms of a digraph on the set of orbits.

\subsubsection{Shapes}

Let $\cF$ be a fusion law.  A \emph{shape} on an axet $(G, X, \tau)$ is an assignment of a $2$-generated $\cF$-axial algebra $A_Y$ to each $2$-generated subaxet $Y$ of $X$, so that $Y$ is naturally identified with $X(A_Y)$.  This assignment should satisfy the obvious consistency requirements.  If we have inclusion between two subaxets $Y \subseteq Z$, then we should have an inclusion of $A_Y \subseteq A_Z$.  Also, if $Y$ and $Z$ are two conjugate subaxets of $X$, where $Y^g = Z$, then we should have a suitable isomorphism $\theta_g \colon A_Y \to A_Z$.  The technical details of this are given in \cite{forbidden}, but to give the reader an idea: a shape is roughly an amalgam of small ($2$-generated) algebras organised on an axet.

Just like every axial algebra contains an axet, every axial algebra includes a shape on its axet.  Thus shapes classify different axial algebras.  Note that shapes rely on us knowing every $2$-generated $\cF$-axial algebra and this is why the investigation of a class of axial algebras starts by classifying the $2$-generated examples.

As mentioned in the beginning, not every shape comes from an $\cF$-axial algebra.  In this case, we say that the shape \emph{collapses}.  This can be made more precise with the definition of a universal algebra \cite{hrs1, r, generic}.  In essence, for a fusion law $\cF$ and a number of generators $k$, there exists a universal $k$-generated $\cF$-axial algebra $A$, which is defined not over a field, but a certain commutative ring.  The field versions of the algebra are obtained by factoring the coefficient ring by a prime ideal.  Therefore, the spectrum of the coefficient ring gives us a variety of all $k$-generated $\cF$-axial algebras.  The shape is an invariant on this variety.  For the fusion law $\cM(\frac{1}{4}, \frac{1}{32})$, there are only eight $2$-generated algebras, the Norton-Sakuma algebras, and so the shape is a discrete invariant and it is constant on every component.

Initial calculations with small groups and the Monster fusion law $\cM(\frac{1}{4}, \frac{1}{32})$, would typically output a single algebra for each shape, with a small number of shapes collapsing.  This suggested that the variety was perhaps $0$-dimensional.  However, Whybrow \cite{whybrow family} gave an example of a $1$-parameter family of $3$-generated $\cM(\frac{1}{4}, \frac{1}{32})$-axial algebras all with the same shape.  This means that components of positive dimensions exist.  This shifts the focus from individual algebras to components which are parametrised families of algebras.

While the initial calculations produced almost no collapsing shapes, subsequent systematic work of Khasraw \cite{kPhD} and Khasraw, M\textsuperscript{c}Inroy and Shpectorov \cite{4trans}, seems to imply that the majority of the shapes for $\cM(\frac{1}{4}, \frac{1}{32})$ collapse.  In \cite{forbidden}, for a family of shapes with the general $\cM(\al, \bt)$ fusion law it is in fact shown that the non-collapsing shapes are a rare exception.

\begin{problem}
Find conditions to impose on shapes to eliminate the majority of collapsing shapes while preserving shapes leading to non-trivial algebras.
\end{problem}

\subsection{Algorithms}

Seress developed an algorithm \cite{seress} and a GAP program for computing algebras using the `resurrection' trick \cite{ipss}. This program was very successful for small groups as long as the algebra was $2$-closed (spanned by axes and products of two axes).  Later the same algorithm was redeveloped and expanded by Pfeiffer and Whybrow \cite{Maddyconstruction}, with the $2$-closed restriction partly removed, and it is implemented in GAP \cite{Maddycode}.  This was used, for example, in the project of Mamontov, Staroletov and Whybrow \cite{msw} to classify minimal $3$-generated axial algebras of Monster type $(\frac{1}{4}, \frac{1}{32})$.

An alternative approach based on the idea of expansions, which overcomes any $m$-closed restriction, was originally explored in GAP by Shpectorov and then properly developed by Rehren.  The algorithm was developed much further for a general fusion law by M\textsuperscript{c}Inroy and Shpectorov in \cite{construction} and a {\sc magma} implementation by M\textsuperscript{c}Inroy is available in \cite{ParAxlAlg}.  This program was able to find much larger algebras and tens of thousands of shapes have been analysed with several hundred non-trivial algebras found.

The main computational bottleneck is working with large-dimensional $G$-modules.  This can be tackled in two ways.  Firstly, by working with the modules decomposed into homogeneous components and preserving these when taking tensor and symmetric products; this means that taking submodules becomes a few much smaller linear algebra problems rather than a large $G$-module problem.  Secondly, by doing the entire calculation over several finite fields and using a Chinese-remainder-theorem-like argument to recover the result for characteristic $0$.  Implemented either or both of these is an interesting computational problem and would allow us to explore much larger algebras.

\begin{problem}
Can these ideas lead to an explicit computed implementation of the Griess algebra?
\end{problem}

\section{Further areas of research}
In this section, we briefly discuss some additional research in this area.

\subsection{Algebras for simple groups}

One of the motivations for group theorists to consider axial algebras is to have a unified approach to simple groups that arise from these algebras.  Classical groups and $G_2$ are realised as automorphism groups of Jordan algebras, whereas the Monster and many of its subgroups act on the Griess algebra and its subalgebras.  So, within the class of axial algebras of Monster type, we have algebras for many simple groups.  It would be very interesting to find axial realisations for the remaining simple groups.

Norton \cite{nortonPhD} proposed a class of non-associative algebras coming from $3$-transposition groups, which are similar to Matsuo algebras.  In \cite{cameron norton alg}, Cameron, Goethals and Seidel reinterpreted Norton's examples in terms of the association scheme coming from a permutation action of a group $G$.  This class of algebras are now called Norton algebras.  Note that Norton algebras are commutative non-associative algebras and they are typically realised from an irreducible component of the permutation module of $G$.

If one wants to have an identity in the algebra, then one can extend the Norton algebra using a cocycle (a $G$-invariant bilinear form).  For example, Norton himself noticed that the Griess algebra can be obtained in this way.  The Norton construction spurred interest in algebras for simple groups, in particular for sporadic groups.  Algebras were constructed for $3.F'_{24}$ by Norton \cite{nortonPhD}, $A_6$, $A_7$, $\Omega^-_6(3)$, $\Omega_7(3)$ and $Fi_{24}'$ by Smith \cite{smith} and Kitazume \cite{kitazume}, $J_3$ \cite{frohardt J3} and $ON$ \cite{frohardt ON} by Frohardt, and $HN$ \cite{Ryba HN} and $B$ \cite{Ryba B} in positive characteristic by Ryba.  However, this was before the axial language was introduced and so the axial structure for these algebras has not been investigated.

\begin{problem}
Explore the above algebras from the axial point of view and construct interesting axial algebras for other groups.
\end{problem}

There has been some progress in this direction.  Recently, Shumba in \cite{ShumbaPhD}, investigated the axial structure of two such algebras of dimensions $1+77$ and $1+780$ for the sporadic groups $HS$ and $Suz$, respectively.  He found the fusion law for specific idempotents and showed that these algebras are axial algebras whose Miyamoto group recovers the original sporadic group. (Strictly speaking these are axial decomposition algebras \cite{dpsv} as to have a grading, we need to split one of the eigenspaces into two separate parts.)  M\textsuperscript{c}Inroy, Peacock and Van Couwenberghe have investigated the axial structure for $J_3$, which has nilpotent `axes' and has a $C_3 \times C_3$-grading, and for $Ly$, which is suspected to have a $C_3$-grading.

In \cite{MichielPhD} and \cite{E8 alg}, De Medts and Van Couwenberghe gave a general procedure for constructing a commutative non-associative algebra for any simply-laced Chevalley group, by looking at a constituent of the symmetric square of the Lie algebra.  They considered these from an axial point of view, calculating the fusion law in all cases and showing that they are axial decomposition algebras.  They concentrated, in particular on the example for $E_8$ of dimension $1+3875$.  In this case, the algebra had been constructed previously from a different point of view by Garibaldi and Guralnick \cite{gar gur}, who were interested in groups stabilising a polynomial on a vector space.

One deficiency in all the algebras in this section, is that the fusion laws tend not to be nice for algebras constructed in this way.  They typically have more eigenvalues than for the Jordan, or Monster fusion laws and some eigenspaces split.  It should be noted, however, that they are typically constructed from a single irreducible module and we can instead consider algebras built in a similar way from a sum of several irreducibles.  In fact, many of the known algebras of Monster type are far from being irreducible.

\begin{question}
Develop methods for constructing algebras on reducible modules.  How can a prescribed fusion law be used to restrict possibilities for the product?
\end{question}


\subsection{Automorphism groups and finding all axes}\label{sec:aut}

From the axial point of view, it is interesting to find all axes in an axial algebra.  On the other hand, from the group-theoretic point of view, it is interesting to find the full automorphism group of an algebra.  It turns out that these questions are intimately related for axial algebras.  Indeed, the automorphism group has a natural permutation action on the full set of axes and conversely the Miyamoto group of the full set of axes is a normal subgroup in the automorphism group.

For simple Jordan algebras the full automorphism groups are the classical groups (and the group $G_2$) and it is known that the Monster is the full automorphism group of the Griess algebra and that the $2\A$-axes are the full set of $\cM(\frac{1}{4}, \frac{1}{32})$-axes.  Castillo-Ramirez found all the idempotents and calculated the automorphism groups for the Norton-Sakuma algebras \cite{NS idempotents} and also for two algebras, of dimensions $6$ and $9$, with the Miyamoto group $S_4$ in \cite{assoc subalgs}.  However, in general, the number of all idempotents in $A$ is either infinite, or finite but growing exponentially, being generically $2^{\dim(A)}$.  Currently there is a project underway with M\textsuperscript{c}Inroy, Shpectorov and Shumba to find computationally the automorphism groups of the remaining algebras for $S_4$, which have dimension up to $25$.  The method involves reducing the problem from $A$ to the smaller subalgebra $A_0(a)$ for an axis $a$.

\begin{problem}
Develop computational methods for finding axes and the full automorphism groups of axial algebras.
\end{problem}

Combining such methods with group-theoretic approaches, we can potentially do much larger algebras -- the Griess algebra was completed using group theory.

Another reason why finding the full automorphism group  is important is that in the list of known examples of algebras, for example in \cite[Tables 4 and 5]{construction} and \cite[Table 4]{Maddyconstruction}, the dimensions tend to repeat.  This suggests that there may be several seemingly different algebras which are, in fact, the same.

Suppose that two axial algebras $A = (A, X)$ and $B = (B,Y)$ are isomorphic via an isomorphism $\phi \colon A \to B$.  Then in $B$, in addition to the closed set of axes (axet) $Y$, there is  also the closed set of axes $\phi(X)$.  Either one set is contained in the other, which is typically easy to check, or neither of them is the full set of axes.  In this case, $B$ has a larger set of axes and perhaps a larger automorphism group.  In most cases, there are not additional axes, however recently Alharbi\footnote{Alharbi is currently a PhD student at the University of Birmingham and this is part of his PhD project.} found three series of Matsuo algebras, with $\eta \neq \frac{1}{2}$, which contain additional axes and have a larger automorphism group (also being a $3$-transposition group).

The above discussion of possibly isomorphic algebras $A$ and $B$ equally applies where the axes from $X$ and $Y$ obey different fusion laws.  This situation does indeed happen.  A simple example is $3\C(\al)$, $\al \neq -1$.  In addition to the three axes $a, b, c$ of Jordan type $\al$, we have three axes $\1-a, \1-b, \1-c$ of Jordan type $1-\al$.  The following example noticed by M\textsuperscript{c}Inroy extends this idea.

\begin{example}
Let $A = 4\B(\al, \frac{\al^2}{2})$ with axes $a_0, \dots, a_3$.  It can be seen that $\lla a_0, a_2 \rra \cong 3\C(\al) \cong \lla a_1, a_3 \rra$.  When $\al \neq -1$, let $\1_+$ be the identity in $\lla a_0, a_2 \rra$ and $\1_-$ be the identity in $\lla a_1, a_3 \rra$. Define $a_i' := \1_+ - a_i$, if $i = 0,2$, and $a_i' := \1_- - a_i$, if $i = 1,3$.  Then the idempotents $a_i'$ are in fact axes of Monster type $(1-\al, \frac{\al(2-\al)}{2})$.  Setting $\gamma = 1-\al$, we check that $ 4\B(\al, \frac{\al^2}{2}) = \lla a_0', a_1' \rra \cong 4\Y(\gamma, \frac{1-\gamma^2}{2})$.  So the two algebras $4\B(\al, \frac{\al^2}{2})$ and $4\Y(\gamma, \frac{1-\gamma^2}{2})$ are isomorphic as algebras, but have different fusion laws and disjoint sets of axes.
\end{example}


\subsection{Code algebras}

Code algebras were introduced by Castillo-Ramirez, M\textsuperscript{c}Inroy and Rehren in \cite{codealgebras}.  They axiomatise some properties found in code VOAs in a similar way that axial algebras do for OZ-type VOAs.  Code VOAs are an important class of VOAs, where a binary code controls the representation theory.

\begin{definition}
Let $C \subseteq \FF_2^n$ be a binary linear code of length $n$ and $\FF$ a field, where $a,b,c \in \FF$.  Set $C^* := C - \{ {\bf 0}, \1\}$.  The \emph{code algebra} $A_C$ is the algebra on $\bigoplus_{i = 1}^n \FF t_i \oplus \bigoplus_{\alpha \in C^*} e^\alpha$ with multiplication given by
\begin{align*}
t_i t_j &= \delta_{i,j} t_i \\
t_i e^\alpha & = \begin{cases} 
a \, e^\alpha & \text{if } \alpha_i = 1 \\
\mathrlap0 \phantom{\sum_{i \in \supp(\alpha)}c\, t_i} & \text{if } \alpha_i =0
\end{cases} \\
e^\alpha e^\beta & = \begin{cases}
b\, e^{\alpha + \beta} & \text{if } \alpha \neq \beta, \beta^c \\
\sum_{i \in \supp(\alpha)}c\, t_i & \text{if } \alpha = \beta  \\
0 & \text{if } \alpha = \beta^c
\end{cases}
\end{align*}
\end{definition}


The $t_i$ are axes of Jordan type $b$, however they do not generate the algebra.  In \cite{codealgebras, CM19} Castillo-Ramirez, M\textsuperscript{c}Inroy and Rehren give a construction, called the $s$-map, which gives an idempotent $s(D,v)$ from a subcode $D \subseteq C$ and a vector $v \in \FF_2^n$.  The axial properties of this idempotent can be difficult to analyse, but this can be done for the smallest such subcodes $\la \alpha \ra$ and we call such idempotent \emph{small idempotents}.  The small idempotents whose fusion law is $C_2$-graded are classified in \cite{CM19} and the resulting Miyamoto groups of the algebras are given in \cite{Miy code alg}.  However, it turns out that these Miyamoto groups for small idempotents are all abelian.

The following example, shows that we do indeed get interesting groups for idempotents coming from the $s$-map construction for larger subcodes $D \subseteq C$.

\begin{example}
Let $C = H_8$ be the Hamming code of length $8$.  Then $A_{H_8}$ is a $22$-dimensional code algebra.  We consider the $s$-map idempotents $s(C,v)$.  Choosing $v \in \FF_2^8 - C$ to have odd weight gives a set of eight mutually orthogonal axes of Jordan type $b$ and choosing $v$ to have even weight gives a different set of eight mutually orthogonal axes of Jordan type $b$.  Both of these are different from the set $\{ t_i : i = 1, \dots, 8\}$ of axes of Jordan type $b$.  These 24 axes together generate the algebra and hence $A_{H_8}$ is an axial algebra of Jordan type $b$ and its Miyamoto group has shape $2^6{:}S_3$.  Moreover, $\Aut(H_8)$ also acts on the algebra, so we see a group of automorphism of shape $2^6{:}(L_3(2) \times S_3)$.
\end{example}

For code VOAs, the example $V_{H_8}$ coming from the Hamming code is very important.  Miyamoto used $V_{H_8}$ in a new construction of the Moonshine VOA $V^\natural$ and other VOAs \cite{VH8}.  Using the parameters $(a,b,c) = (\frac{1}{4}, \frac{1}{2}, 1)$ above, $A_C$ is isomorphic to the weight $2$ part of $V_{H_8}$ and moreover, we see the full automorphism group $\Aut(V_{H_8}) = 2^6{:}(L_3(2) \times S_3)$ as a group of automorphisms of $A_C$.

\end{document}